\documentstyle[11pt]{article}

\input amssym.def

\newtheorem{lemma}{Lemma}[section] 
\newtheorem{proposition}[lemma]{Proposition} 
 
\newtheorem{theorem}[lemma]{Theorem} 

\newtheorem{defn}[lemma]{Definition}
\newtheorem{exam}{Example}

\newtheorem{rem}{Remark}

\def\Sing{\mathop{\rm Sing}\nolimits} 
\def\Aut{\mathop{\rm Aut}\nolimits} 
\def\Map{\mathop{\rm Map}\nolimits} 
 
\def\ev{\mathop{\rm ev}\nolimits}

\def\Euler{\mathop{\rm Euler}\nolimits} 
\def\spin{\mathop{\rm spin}\nolimits} 

\def\C{{\cal C}}

\def\E{{\cal E}}
\def\M{{\cal M}}
\def\F{{\cal F}}
\def\B{{\cal B}}
\def\U{{\cal U}}
\def\B{{\cal B}}
\def\A{{\cal A}}
\def\G{{\cal G}}
\def\L{{\cal L}}
\def\F{{\cal F}}

\newcommand{\CC}{{\Bbb C}}
\newcommand{\QQ}{{\Bbb Q}}
\newcommand{\ZZ}{{\Bbb Z}}
\newcommand{\RR}{{\Bbb R}}

\def\sl{\frak{sl}}

\def\wp{W^+}
\def\wm{W^-}
\def\Lamp{\Lambda^+}
\let\lra=\longrightarrow

\def\uso{^{S^1}}
\def\bso{ \B^{\so} } 
\def\so{S^1}
\let\sub=\subset

\title{Virtual moduli cycles and Gromov-Witten\\ 
invariants of general symplectic manifolds}

\author{Jun Li\\Department of Mathematics\\
Stanford University\\
Palo Alto, CA 94305 \\ and \\
Gang Tian\thanks{Supported partially by 
NSF grants}\\Department of Mathematics
\\Massachusetts Institute of Technology\\
Cambridge, MA 02139}
\date{}

\begin{document}

%%%%%%%%%%%%%%%%%%%%%%%%%%%%%%%%%%%%%%%%%%%%%%%
%%%%%% The following needs to be changed for latex 2.09
%%%%%%%%%%%%%%%%%%%%%%%%%%%%%%%%%%%%%%%%%%%%%%%
%%%\title{}
%%%%    Information for first author
%%%\author{Gang Tian}
%%%%    Address of record for the research reported here
%%%\address{Massachusetts Institute of Technology, Cambridge, MA
%%%  02139}
%%%%    Current address
%%%\curraddr{Department of Mathematics, Massachusetts Institute of
%%%  Technology, Cambridge, Mass.\ 02139}
%%%
%%%\email{tian@math.mit.edu}

%\section{Introduction}
%\author{Gang Tian \\
%{\small Massachusetts Institute of Technology, Cambridge, MA~
%% 02139} \\
%{\small tian@math.mit.edu}
%}

%\thanks{}

%%%\author{}
%%%\address{}
%%%\email{}

%    General info
%%%\subjclass{Subject classes here}

%%%\date{}

%\dedicatory{Do you want a dedicatory here?}

\maketitle
\begin{abstract}
In this paper, we first give an intersection theory for moduli problems
for nonlinear elliptic operators with certain precompact
space of solutions in differential geometry. Then we apply the theory to
constructing Gromov-Witten invariants for general symplectic manifolds.
\end{abstract}

\setcounter{section}{-1}

\section{Introduction}
\label{sec:0}

This paper is a sequel to [LT1]. In [LT1], by
using purely algebraic methods, we developped an intersection
theory for moduli problems on smooth algebraic varieties
over any algebraically closed fields of characteristic zero. 
An alternative construction was given in [BF].
The key point in [LT1] is the existence of locally 
free resolutions of tangent complexes involved.
In this paper, we apply the same idea to constructing 
the intersection
theory for moduli problems in the differential category.
However, the tool will be analytic this time.

Given a Banach manifold $B$, a smooth bundle $E\mapsto B$ is 
Fredholm, if there is a section $s: B \mapsto E$ such that $s^{-1}(0)$
is compact and any linearization of $s$ at any point of 
$s^{-1}(0)$ is Fredholm
of index $r$. Then one can define the determinant line bundle $\det(s)$
over $s^{-1}(0)$. Note that $s:B\mapsto E$ is orientable if $\det (s)$ does.
It should be well-known that for any orientable Fredholm
bundle $s: B \mapsto E$, one can associate an Euler
class $e([s:B\mapsto E])$ in $H_r(B, \ZZ )$, 
which depends only on the homotopy 
class $[s:B\mapsto E]$ of $s: B\mapsto E$.

However, its use is very limited. In many problems, such as constructing
Donaldson invariants and Gromov-Witten invariants,
the zero set $s^{-1}(0)$ is often noncompact, if we insist on
smooth Banach manifolds, smooth bundles. For 
many useful applications, we have to construct 
Euler classes for spaces, bundles and their sections, which are 
not necessarily smooth, and prove that they are invariant. 
In section one, we will
give two simple theorems on constructing the Euler classes of
so called generalized Fredholm bundles, and more generally, 
the rational Euler classes of generalized
Fredholm orbifold bundles.
The main part of this paper is devoted to constructing
Gromov-Witten invariants over rational numbers 
for general symplectic manifolds by 
establishing the Fredholm properties of the bundle of $(0,1)$-forms
over the space of smooth stable maps (cf. section 2).
In fact, we have constructed symplectically invariant Euler classes
in the space of stable maps.

The theory of the Gromov-Witten invariants
was first established in a systematical and mathematical way by
Ruan and the second author in [RT1], [RT2] for semi-positive
symplectic manifolds. They actually constructed the invariants
over integers. 
In [LT1] and [BF], the authors constructed
the Gromov-Witten invariants for any algebraic manifolds
over any closed fields of characteristic zero.

The similar idea was also used by Liu and the second author in
proving the Arnold conjecture for nondegenerate Hamiltonian
functions on general symplectic manifolds [LiuT].

During the preparation of this paper, we learned that
Fukaya and Ono also gave a different construction of
the Gromov-Witten invariants and a proof of the
Arnold conjecture for nondegenerate Hamiltonian functions
for general symplectic manifolds ([FO]). We also learned
that one or both of them was also claimed
by Hofer and Salamon, and Ruan.
Shortly after we finished writting of this paper, we received
a preprint [Si] of B. Siebert, in which he gave another different
construction of Gromov-Witten invariants for general symplectic
manifolds. 

We believe that our construction can be also used to constructing
the Gauge theory invariants, such as Donaldson invariants. 
It is also interesting to compare 
the Gromov-Witten invariants constructed here with algebraic ones in
[LT1] (cf. [LT2]). We plan to discuss these in forthcoming papers. 
 
\tableofcontents

%%%\clearpage

\section{Euler classes for Fredholm bundles}

In this section, we collect a few simple theorems, which can be proved
easily.

Let B be a topological space. Recall that 
a topological bundle $\pi : E \mapsto B$ consists of a continuous
map $\pi$ between topological spaces, satisfying:
(1) there is a topological subspace $Z$ in $E$ such that $\pi|_Z$ is
a homeomorphism from $Z$ onto $B$; (2) For any $x\in B$, the fiber
$E_x= \pi^{-1}(x)$ is a vector space with origin at $(\pi|_Z)^{-1}(x)$.
A section of $E$ is a map $s: B \mapsto E$ such that
$\pi \cdot s$ is the identity map of $B$.
Clearly, $Z$ defines a section of $E$, which is 
usually refered as the 0-section.
For any section $s$, its zero locus in $B$ is
$s^{-1}(Z)$, which is also denoted by $s^{-1}(0)$.

A smooth approximation $(E_i, U_i)$ of
$s: B\mapsto E$ consists of an open subset $U_i$ in $B$ and a
continuous vector subbundle $E_i$ of finite rank over $U_i$,
such that $s^{-1}(E_i)\subset U_i$ is a smooth manifold and
$E_i | _{ s^{-1}(E_i)}$ is a smooth bundle over 
$ s^{-1}(E_i)$ with $s|_{s^{-1}(E_i)}$ smooth .
We say that $s: B \mapsto E$ has a weakly smooth structure
$\{(E_i, U_i)\}$ of index $r$, if (i) each $(E_i, U_i)$ is a smooth
approximation of $s:B \mapsto E$; 
(ii) $\{U_i\}$ is a covering of $ s^{-1}(0)$;
(iii) $s^{-1}(E_i)\subset U_i$ is of dimension $r + {\rm rk }(E_i)$; 
(iv) For any $x\in s^{-1}(0)\cap 
U_i\cap U_j$, there is another smooth approximation $E_k \mapsto
U_k$ with $x \in U_k$,
such that $E_i|_{U_i\cap U_k}$ (resp. $E_j|_{U_j\cap U_k}$)
is a subbundle of $E_k|_{U_i\cap U_k}$ 
(resp. $E_k|_{U_j\cap U_k}$).

Given two smooth structures $\{(E_i,U_i)\}$ and $\{(E_j', U_j')\}$ of 
$s: B\mapsto E$,
we say that $\{U_i\}$ is finer than $\{U_j'\}$, if
for $x\in s^{-1}(0)\cap U'_j$, there is at least one $(E_i,U_i)$
such that near $x$, (1) $s^{-1}(E_i)\cap s^{-1}(E_j')$ is a smooth 
submanifold in $U_i$ of dimension
$\dim s^{-1}(E_j')$; (2)
$E_j'|_{U_i\cap U_j'}$ is a 
subbundle of $E_i|_{U_i\cap U_j'}$; (3) The restriction
$E_j'|_{s^{-1}(E_i)\cap s^{-1}(E_j')}$ is a smooth subbundle
of $E_i|_{s^{-1}(E_i)\cap s^{-1}(E_j')}$.

We say that $E$ is a generalized Fredholm bundle of index $r$,
if there is a continuous section
$s: B \mapsto E $ satisfying the followings: 

\noindent
(1) $s^{-1}(0)$ is compact; 

\noindent
(2) $s: B\mapsto E$ has a weakly smooth structure $\{(E_i, U_i)\}$
of index $r$;

\noindent
(3) There is a finitely dimensional vector
space $F$ and a bundle homomorphism
$\psi _F: B\times F \mapsto E$, 
such that for any smooth approximation $(E_i, U_i)$, 
$\psi_F|_{s^{-1}(E_i)\times F}$ is a smooth map from
$s^{-1}(E_i)\times F$ into
$E_i$ and transverse to $s$ along $s^{-1}(0)\cap U_i$. 

Such a section $s$ is called admissible. We call $\{F, (E_i, U_i)\}$ a 
weakly smooth resolution of $s: B\mapsto E$.

Put $W_F = (s-\psi_F)^{-1}(0)\subset B\times F$.
Here by abusing the notation, we also regard $s$ as
a section of $E$ over $B\times F$.
Then $W_F$ is a smooth manifold of
dimension $r + {\rm rk}(F)$ near $s^{-1}(0)$, 
and $s$ lifts to a smooth section $s_F: W_F \mapsto
F$, namely, for any $(x, v)\in W_F\subset B\times F$,
$s_F(x,v)= v $. Clearly, $s_F^{-1}(0)= s^{-1}(0)$.

\begin{rem}

One can also define the weakly $C^\ell$-smoothness of
$s:B\mapsto E$. We say that $s:B \mapsto E$ is of class
$C^\ell$ ($\ell \ge 1$) if any $E_i$ is a 
$C^\ell$-smooth bundle over a $C^\ell$-smooth manifold
$s^{-1}(E_i)$, and $s$, $\psi_F$ are $C^\ell$-smooth 
along $s^{-1}(E_i)$.

\end{rem}

\begin{rem}

If $\{F', (E_j', U_j')\}$ is another smooth resolution of $s$, 
we say that $\{F, (E_i, U_i)\}$ is finer than $\{F', (E_j', U_j')\}$,
if $F' \subset F$, $\{(E_i, U_i)\}$ is finer than $\{(E_j', U_j')\}$
and $\psi_{F}$ restricts to $\psi_{F'}$ on $F'$.

We will identify $\{F',(E_j', U_j')\}$ with 
$\{F, (E_i, U_i)\}$, if there is another smooth structure
$\{F'', (E_k'', U_k'')\}$ finer than 
$\{F',(E_j', U_j')\}$ and $\{F, (E_i, U_i)\}$.

Let $s, s' : B \mapsto E$ be two admissible sections.
In the following, unless specified, by $s=s'$, we mean
that they are the same as continuous sections
and their weakly smooth structures are identical.

\end{rem}

We say that two generalized Fredholm bundles
$s: B \mapsto E$ and $ s': B \mapsto E$ 
are homotopic to each other, if there is a generalized
Fredholm bundle of the form $S: \pi^*_2 E \to [0,1] \times B$
of index $r+1$, such that $S |_{0
\times B} = s$ and $S |_{1 \times B} = s'$,
where $\pi_2: [0,1]\times B \mapsto B$ is the natural projection. 

We denote by $[s : B \to E]$ the equivalence class of generalized
Fredholm bundles which are homotopic to $s : B \to E$. 
We also denote by $r(B,E,s)$ the index
of the generalized Fredholm bundle $s: B \mapsto E$. 

\begin{exam}
Let $B$ be a Banach manifold (possibly incomplete). 
Suppose that $E$ is a vector bundle $E$ over $B$ with 
a section $s : B \to E$ satisfying: 
1) $s^{-1} (0)$ is compact; 2) For any $x \in B$, $L_x
(s)$ is Fredholm, where $L_x (s)$ denotes the linearization of
$s$ at $x$ with $s(x) =0$, defined as follows:
if $\phi : E|_W\mapsto W\times H$ is any local trivialization near $x$, 
then 
$$L_x(s) (v) = \phi^{-1} (v (\pi_2\cdot \phi \cdot s)(x)),$$
where $\pi_2$ be the projection from
$W\times H$ onto $H$. Since $s(x) =0$, $L_x(s)(v)$ is independent
of choices of local trivializations of $E$ near $x$.

The index of $L_x (s)$ is independent of $x$ in $B$. 
Therefore, we can denote
this index by $r (B, E, s)$.

One can easily show that such a $s: B \mapsto E$ is a generalized Fredholm
bundle of index $r(B, E, s)$.

\end{exam}

Let $s: B\mapsto E$ be a generalized Fredholm bundle.
We can define its determinant bundle $\det (s)$ as follows:
let $\{F, (E_i, U_i)\}$ be a smooth resolution of $s$, then we define
$\det (s)$ to be $\det(TW_F)\otimes \det(F)^{-1}|_{s^{-1}(0)}$.
For any smooth approximation $(E_i, U_i)$, $W_F $ is 
a smooth submanifold in $ U_i \times F$ and its normal
bundle can be canonically identified with $E_i|_{W_F}$
by using the differential $d (s- \psi_F)$. It follows that
$det(s) $ can be canonically identified with
$\det (TU_i\times F)\otimes \det (E_i)^{-1} \otimes \det (F)^{-1}|_{s^{-1}(0)}$, 
and consequently, $\det (TU_i)\otimes \det (E_i)^{-1}|_{s^{-1}(0)}$.
In particular, $\det(s)$ is independent of choices of $\{F, (E_i, U_i)\}$.
Moreover, it implies that $\det (s)(x)$ ($x\in s^{-1}(0)$) can be
naturally identified with $\det (L_x(s))$, where
$L_x(s)$ denotes the linearization of $s$ from $\bigcup \{T_xU_i |
x \in U_i\cap s^{-1}(0)\}$ into $E_x$ defined as in Example 1.
The determinant $\det (L_x(s))$ is defined in the standard way
by using finitely dimensional approximations.

We say that $s: B \mapsto E$ is orientable if 
${\rm det}(s)$ is orientable, i.e., it admits a
nonvanishing section. Clearly, if $s: B\mapsto E$ is orientable,
so is any other bundle in $[s:B\mapsto E]$, so we can simply say that
$[s:B\mapsto E]$ is an orientable equivalence class.

\begin{theorem}

  For each oriented equivalence class $[s: B\to E]$ of
  generalized Fredholm bundles, we can assign an
  oriented Euler class $e([s:B\to E])$ in 
  $ H_r(B, Z)$, where $r=r(B,E,s)$. More 
  precisely, $e([s : B \to E])$ can be represented by 
  an $r$-dimensional manifold.

  Furthermore, this Euler class satisfies the usual functorial properties
  for the Euler class of bundles of finite rank.

\end{theorem}
\vskip 0.1in
\noindent
{\bf Proof:} 
First we observe that $s_F^{-1}(0) = s^{-1}(0)$ is compact.
Then, by the standard transversality theorem, there is a generic,
small section $v: W_F \mapsto F$, such that 
$(s_F+v)^{-1}(0)$ is a compact submanifold in $B\times F$ of dimension $r$. 
It has a natural orientation induced by $\det(s)$.

We claim that the homology class of $(s_F + v)^{-1}(0)$ is
independent of choices of smooth resolutions.
Suppose that $\{F', (E_j', U_j')\}$ is another smooth
resolution which is finer than $\{F, (E_i, U_i)\}$.
Then we have another smooth manifold $W_{F'}$ containing
$W_F$ as a submanifold. 

We may assume that 
$W_F = W_{F'} \cap B\times F$ and the above $v$ extends
to a map from $W_{F'}$ into $F$.

Let $F' = F\oplus F'/F$ be a splitting. 
Write $B\times F'$ as $B\times F\times F'/F$, 
for any $(x,u_1, u_2)\in W_{F'}\backslash W_{F}$, we have
$u_2 \not= 0$. It follows that $(s_{F'} + v)^{-1}(0) = (s_F+v)^{-1}(0)$. 
Then the claim follows easily.

We define $e([s: B\mapsto E])$ to be the homology class
in $H_r(B, \ZZ)$, which is  
represented by $(s_F + v)^{-1}(0)$ in $B\times F$.
Here we have used the fact that $B$ is homotopically equivalent
to $B\times F$.
We can also regard $e([s:B \mapsto E])$ as the intersection
class of $W_F$ with $B\times \{0\}$.

The class $e([s:B \mapsto E])$ is independent 
of choices of Fredholm sections $s$
in $[s: B \mapsto E]$. In fact, to prove it, we simply
repeat the above arguments for any homotopy
$S: B\times [0,1] \mapsto E$ with $S|_0 = s$.

One can show that $e([s: B \mapsto E])$ satisfies all
functorial properties of the Euler class. So Theorem 1.1 is proved.

\begin{rem}

Theorem 1.1 still holds even if the assumption (3) on $s: B\mapsto E$
is replaced by

\noindent
(3)' there are finitely many 
open subsets $\{V_a\}$ and finitely dimensional vector spaces
$F_a$ satisfying: (i) $B = \bigcup _a V_a$; (ii)
For each $a$, there is a bundle map $\psi_a: V_a\times F_a \mapsto E|_{V_a}$
which is transverse to $s$ along $s^{-1}(0)$ for
any smooth approximations, or equivalently, $(F_a, \psi_a)$ is
a weakly smooth resolution of $s|_{V_a}$;
(iii) If $\dim F_a \le \dim F_b$, 
we have that $F_a \subset F_b$ and 
$$\psi_b |_{(V_a \cap V_b )\times F_a} = \psi _a |_{(V_a \cap V_b )\times F_a}.$$
The proof is not much more difficult than the above one.

We can also regard $\{V_a, F_a\}$ as a resolution of $s$. However, when
the weakly smooth structure of $B$ admits appropriate cut-off functions, 
this weakened condition is the same as (3).
\end{rem}

The assumptions in Theorem 1.1 can be weaken, namely, we do not really
need $B$ to be weakly smooth. The following is rather straightforward, if
one treats orbifolds like manifolds. 

Let $B$ be a topological space.
We recall that a topological fibration $\pi: E\mapsto B$
is an orbifold bundle, if there is a covering $\{ U_i\}$ of $B$ by open subsets,
satisfying: (1) each $U_i$ is of the form $\tilde U_i / \Gamma_i$, where
$\Gamma _i$ is a finite group acting on $\tilde U_i$;
(2) for each $i$, there is a topological bundle
$\tilde E_i \mapsto \tilde U_i$, such that $E|_{U_i}=\tilde E_i/ \Gamma_i$;
(3) For any $i, j$, there is a bundle map $\Phi_{ij}$
from $\tilde E_j|_{\pi^{-1}_j(U_i\cap U_j)}$ to 
$\tilde E_i|_{\pi^{-1}_i(U_i\cap U_j)}$, which descends to
the identity map of 
$E|_{U_i \cap U_j}$,
where $\pi_k : \tilde U_k \mapsto U_k$ is the natural projection;
(4) For each $x\in \pi_j^{-1}(U_i\cap U_j)$, there is a
small neighborhood $U_x$, such that 
$\Phi_{ij}|_{\pi_j^{-1}(U_x)}$ is an isomorphism from
each connected component of $\pi_j^{-1}(U_x)$ onto its image;
(5) Each $\Phi_{ij}$ is compatible with actions of
$\Gamma_i$ and $\Gamma_j$.
Any such a $\pi_i : \tilde U_i\mapsto U_i$ is called a
local uniformization of $B$.
Note that $\Phi_{ij}\cdot \Phi_{ji}$ may not be an identity. It is only
a covering map. We will denote by $\phi_{ij}$ the induced map from
$\pi_j^{-1}(U_i\cap U_j)$ to $\pi_i^{-1}(U_i\cap U_j)$.

An orbifold section $s: B\mapsto E$ is a continuous map such that
for each $i$, $s|_{U_i}$ lifts to a section $s_i$
of $\tilde E_i$ over $\tilde U_i$.

Similarly, we can define orbifold bundle homomorphisms,
and the zero set $s^{-1}(0)$ in an obvious way.

Let $\pi: E\mapsto B$ be a topological orbifold bundle.
We say that $E$ is a generalized Fredholm orbifold bundle of index
$r$, if there is an orbifold section $s:B\mapsto E$ satisfying:

\noindent 
(1) $s^{-1}(0)$ is compact; 

\noindent
(2) For each local uniformization $\pi_i: \tilde U_i \mapsto U_i$, 
$s_i$ has a weakly smooth structure of index $r$;

\noindent
(3) For any $i,j$, $\Phi_{ij}$ respects weakly smooth structures
of $s_j: \tilde U_j \mapsto \tilde E_j$ and
$s_i: \tilde U_i \mapsto \tilde E_i$;

\noindent
(4) For each $i$, there is a finitely dimensional vector space
$F_i$, on which $\Gamma _i$ acts, and a $\Gamma_i$-equivariant 
bundle homomorphism $\psi_{i}:\tilde U_i \times F_i\mapsto \tilde E_i$, 
satisfying: (i) together with the weakly smooth structure,
$F_i$ provides a weakly smooth resolution of $s_i$;
(ii)  For any pair $ i, j$, if $\dim F_j \le \dim F_i$,
then there is an injective
bundle homomorphism $\theta _{ij}: \pi_j^{-1} (U_i\cap U_j)\times F_j
\mapsto \pi_i^{-1} (U_i\cap U_j)\times F_i$, such that
$\tilde p _i\cdot \theta_{ij}= \phi_{ij} \cdot \tilde p_j$, where 
$\tilde p_i$ denotes the natural projection from $\tilde U_i \times F_i$ 
onto $\tilde U_i$, and $\psi_i \cdot \theta_{ij} = \Phi_{ij}\cdot
\psi_j $ on $\pi_j^{-1}(U_i\cap U_j)\times F_j$;
(iv) If $\dim F_k \le \dim F_j \le \dim F_i$, then $\theta _{ik}=\theta _{ij}\cdot
\theta_{jk}$ over $\pi_k^{-1}(U_i\cap U_j)$;
(v) For any $x$ in $U_i\cap U_j$, $\theta_{ij}$ is $\Gamma_x$-equivariant
near $\pi^{-1}_j(x)$, where $\Gamma_x$ is the uniformization group of 
$B$ at $x$. We will also call $\{F_i, \psi_i\}$ a resolution of $s$.

Clearly, all $\tilde U_i\times F_i / \Gamma_i$ can be glued together
to obtain a topological space $V(F)$. There is a natural projection
$p_F: V(F)\mapsto B$. In fact, $V(F)$ is a union of finitely many orbifold bundles,
so we may call it an orbifold quasi-bundle. Also, all $\psi_i$ can be 
put together to form a map $\psi_F$ from $V(F) $ into $E$. 

Similarly, one can define notions of homotopy equivalence
of generalized Fredholm orbifold bundles. One can also compare
weakly smooth structures and resolutions of generalized
Fredholm orbifold bundles in the same way as we did before.

For any generalized Fredholm orbifold bundle
$s: B \mapsto E$, we can also associate a determinant orbifold
bundle, denoted by $\det (s)$, in the same way as we did before.
We say that $s: B \mapsto E$ is orientable if
$\det (s) $ does, i.e., $\det (s)$ admits
a nonvanishing section. Note that the orientability of generalized
Fredholm orbifold bundles is a homotopy invariant.

Now we have the following generalization of Theorem 1.1.

\begin{theorem}
\label{theo:1.2}
  For each equivalence class $[s: B\to E]$ of
  generalized Fredholm orbifold bundles, we can assign an
  oriented Euler class $e([s:B\to E])$ in 
  $H_r(B, \QQ )$, where $r$ is the index.

  Furthermore, this Euler class satisfies the usual functorial properties
  for the Euler class of bundles of finite rank.

\end{theorem}

\vskip 0.1in
\noindent
{\bf Proof:} We will adopt the notations in the above definitions
of generalized Fredholm orbifold bundles. 

As before, we denote by $W_F$ the zero set $(s- \psi_F)^{-1}(0)$.
Here again, we regard $s$ as a section of $E$ over $B\times F$
in an obvious way.
For each local uniformization $\pi_i : \tilde U_i \mapsto U_i$,
$\tilde W_{i}= (s_i - \psi _{i})^{-1}(0)$ is smooth
in $\tilde U_i \times F_i$ near $s_i^{-1}(0)$ and
of dimension $r + {\rm rk}(F_i)$. Then 
$W_F$ is obtained by gluing together all $W_i=\tilde W_{i}/ \Gamma_i$. 
More precisely, for any $\dim F_j \le \dim F_i$, by the above (4),
$\tilde W_j \cap \left ( \pi_j^{-1}(U_i) \times F_j\right )$
is locally embedded into $\tilde W_i$ by $\theta_{ij}$,
consequently, we can identify 
$\left (\tilde W_j\cap \pi_j^{-1}(U_i)\times F_j\right )/\Gamma_j$
with a smooth suborbifold, simply denoted by $W_i\cap W_j$ if there is no confusion,
in $W_i$. We also denote by $\pi_i$ the projection from
$\tilde W_i$ onto $W_i$.

Furthermore, $V(F)$ pulls back to an orbifold quasi-bundle $V_0(F)$ over $W_F$.
More precisely, $V_0(F)= \bigcup _i V_i$, where each $V_i= \tilde V_i/\Gamma_i$ 
and $\tilde V_{i} = \tilde W_{i}\times F_i$. As above, 
if $\dim F_j \le \dim F_i$, $\theta _{ij}$ induces an injective bundle map 
$$h_{ij}: V_j|_{W_i\cap W_j} \mapsto V_i|_{W_i\cap W_j}.$$
We denote by $F_{ij}$ the orbifold subbundle
$h_{ij} (V_j|_{W_i\cap W_j})$. It extends to
an orbifold subbundle, still denoted by $F_{ij}$,
of $V_i$ over a small neighborhood
of $W_i\cap W_j$ in $W_i$. Let $\tilde F_{ij}$ be the
lifting of $F_{ij}$ in $\tilde V_i$ over a small neighborhood
of $\pi_i^{-1}(W_i\cap W_j)$ in $\tilde W_i$. 

We define a continuous section 
$s_F: W_F \mapsto V_0(F)$ as follows: for each $i$,
$s_F|_{W_i}$ is descended from the section
$$ (x, v) \in \tilde W_i \subset \tilde U_i \times F_i \mapsto v \in F_i.$$
Clearly, $s_F^{-1}(0) = s^{-1}(0)$.
Moreover, by the definitions of $s_F$ and $W_F$, we may assume that
for some small neighborhood $O$ of $s_F^{-1}(0)$ in $W_F$
and any $\dim F_j \le \dim F_i$,
$F_{ij}$ is well defined over $O\cap W_i\cap p_F^{-1}(U_j)$ and
$$\left ( s_F|_{W_i}\right )^{-1}(F_{ij}|_{O\cap W_i\cap p_F^{-1}(U_j)}) = 
\left (s_F|_{W_j}\right )^{-1} (V_j|_{W_i\cap W_j}).$$

To save the notations, we will simply identify $W_F$ with 
the 0-section in $V_0(F)$.

Let $p: V_0(F) \mapsto W_F$ be the natural projection induced by
$p_F$. Observe that it induces a homomorphism 
$\tau : H_*(V_0(F), \QQ) \mapsto H_*(B,\QQ)$. 
Then we define the Euler class $e([s:B\mapsto E])$ to
be $\tau ([W_F \cap G(s_F)])$. Here,
$[W_F \cap G(s_F)]$ denotes the intersection
class of $W_F$ with the 
graph $G(s_F)$ of $s_F$ in $V_0(F)$.
Using the above properties of $s_F$ and 
standard arguments, one can show that such an intersection exists in
$H_*(V_0(F), \QQ )$. Note that $W_F$ and $V_0(F)$ are unions
of finitely many orbifolds and $s^{-1}_F(0)$ is compact. 

For the reader's convenience, we will outline  
construction of this intersection class by constructing its rational 
cycle representative.

Choose $U_i'$ such that $s^{-1}_F(0)\subset \bigcup _i U_i'$
and its closure $\overline U_i'$ is contained in $U_i$.
Put $\tilde W_i' = \tilde W_i \cap \left (\pi_i^{-1}(U_i')\times F_i\right )$
and $W_i'= \tilde W_i'/ \Gamma_i$. Then
$s_F^{-1}(0) \subset \bigcup _i W_i'= W_F'$.
We also put $\tilde V_i'= \tilde W_i'\times F_i$ and
$V_i'=\tilde V_i'/\Gamma_i$. 

In this proof, by a cocycle $Z'$ of degree $m$ 
in $p^{-1}(O)$, we mean a union
of cycles $Z_i'\subset \overline V_i'\cap p^{-1}(O)$ with its boundary
in $\partial \overline V_i'$ and of dimension $\dim F_i + m$, 
such that $Z_j'\cap \overline V_i'$ is
embedded in $Z_i'$ whenever $\dim F_j \le \dim F_i$,
and $Z_i'\cap p^{-1}(U) \subset F_{ij}$, where $U$ is
some neighborhood of $O\cap \overline W_i' \cap p_F^{-1}(\overline U_j')$.
For example, $W_F$ and $G(s_F)$ are cocycles of degree $r$. 
We say that two cocycles are homologous if they can be
deformed to each other through a family of cocycles.

Let $m_i$ be the order of the local uniformization group
$\Gamma _i$ of $V_i$. Put $m = m_1\cdots m_\ell$.
Then $m e([s:B\mapsto E])$ should be the intersection
class of $m \,W_F$ with the graph $G(s_F)$
in $V_0(F)$. 

We will construct an oriented cocycle $Z$ in $p^{-1}(O)$, which is
homologous to $m W_F$ in $p^{-1}(O)$, such that
$Z$ is transverse to the graph of $s_F$ in $V_0(F)$. 

We will use the induction for this purpose. 

Without loss of generality, we may arrange 
$$\dim F_1 \le \dim F_2 \le \cdots \le \dim F_k \le \cdots.$$

By perturbing $\tilde W_1$ in $\tilde W_1 \times F_1$ and averaging
over the action of $\Gamma_1$, we obtain a cycle 
$\tilde Z_1\subset \tilde V_1=\tilde W_1\times F_1$, satisfying:
(i) $\tilde Z_1$ is $\Gamma_1$-invariant;
(ii) $\tilde Z_1= m_1 \tilde W_1$ near $\partial \tilde V_1$;
(iii) $\tilde Z_1$ is homologous to $m_1 \tilde W_1$ in $\pi_1^{-1}(O)$
with fixed boundary;
(iv) $\tilde Z_1$ is transverse to $G(s_{F,1})$ in an neighborhood
of $\overline {\tilde W_i'\times F_1}$, where 
$s_{F,i}$ be the induced section over $\tilde W_{i}$ by $s_F$.

We extend $\tilde Z_1 /\Gamma_1$ to a cycle over 
$\left (\bigcup_{i\ge 2} O\cap W_i\cap p_F^{-1}(U_1)\right )\cup W_1$, 
such that $\tilde Z_1/\Gamma _1$ is contained in $F_{i1}$
over some neighborhood of $O\cap \overline { W_i'\cap p_F^{-1}(U_1')}$ and 
coincides with $m_1 W_F$ near each 
$\left (\partial (O\cap W_i\cap p_F^{-1}(U_1))\right )\cap W_i$.
Then we can glue $\tilde Z_1/\Gamma_1$ and $m_1 W_F$
together to form a cocycle $Z_1$ in $p^{-1}(O)$.

Now suppose that for $k \ge 1$,
we have found cocycles $Z_i $ ($1\le i \le k$) in $p^{-1}(O)$ formed
by glueing $m_i Z_{i-1}\backslash p^{-1}(W_i)$ and $\tilde Z_i/\Gamma_i$,
where $\tilde Z_i \subset \tilde V_i = \tilde W_i\times F_i$,
satisfying:

\noindent
(i) $\tilde Z_i$ is $\Gamma_i$-invariant
and coincides with $\pi^{-1}_i(Z_{i-1})$, where 
$\pi_i: \tilde V_i \mapsto V_i$ is the projection, near
$\partial \tilde V_i$ and $\pi_i^{-1}(\bigcup_{j< i} p^{-1}( W_j))$;

\noindent
(ii) $\tilde Z_i$ is transverse to 
the graph $G(s_{F,i})$ in an neighborhood of $\overline {\tilde W_i'\times F_i}$;

\noindent
(iii) $\tilde Z_i$ is homologous to $\pi_i^{-1}( Z_{i-1})$
in $\tilde V_i$ with boundary fixed.

Furthermore, we may assume that for each $l > i$,
$Z_i$ is contained in $ F_{li}$
over some neighborhood of $O\cap \overline {W_l'\cap U_i'}$.

Now we construct $Z_{k+1}$.
Observe that $p|_{Z_k}$ is a branched covering of $Z_k$ over $W_F$
of order $m_1\cdots m_k$. It follows that 
$\pi_{k+1}^{-1}(Z_k) \subset \tilde V_{k+1}$, counted with multiplicity,
is a branched covering of $ W_{k+1}\subset W_F$ 
of order $m_1\cdots m_{k+1}$.
Then, by the standard transversality theorem and the same arguments
as we did for $Z_1$,
we can have a $\Gamma_{k+1}$-invariant perturbation $\tilde Z_{k+1}$ of 
$\pi_{k+1}^{-1}(Z_k)$ inside $p^{-1}(O)$,
such that all properties for $\tilde Z_i$ ($i\le k$) are satisfied for
$\tilde Z_{k+1}$.
We define $Z_{k+1}$ to be the glueing of $\tilde Z_{k+1}/ \Gamma_{k+1}$ with
$m_{k+1} Z_k$ along the boundary. The method is the same as that in
the construction of $Z_1$, so we omit it.

The orientation of each $\tilde Z_i$ ($\ge 1$)
is naturally induced by that of $\det (s)$, as we did
before. 

Thus by induction, we have constructed an integral cocycle $Z$ in 
$V_0(F)$ of degree $r$, homologous to $m W_F$ as we wanted.
Moreover, one can show that the intersection of $G(s_F)$ with
$Z$ is a well-defined cycle in $V_0(F)$.
 
We define $e([s: B\mapsto E])$ to be the homology class
in $H_r(B, \QQ )$, which is  
represented by $\frac{1}{m} Z\cap G(s_F)$.
We remind the readers that $G(s_F)$ is the graph of $s_F$ in $V_0(F)$.
This class is independent of choice of the admissible orbifold 
section $s$ in $[s: B \mapsto E]$. In fact, to prove it, we simply
repeat the above arguments for any homotopy
$S: [0,1] \times B \mapsto E$ with $S|_{0\times B} = s$.

One can show that $e([s: B \mapsto E])$ satisfies all
functorial properties of the Euler class. So Theorem 1.2 is proved.
\vskip 0.1in

\begin{rem}

It is very important to know when
$e([s: B \mapsto E])$ is an integer class.

Let us stratify $B$ according to the local unformization group, namely,
write $B$ as a disjoint union of $B_i$, such that the
local unformization group is the same at any points of $B_i$.
In fact, each $B_i$ consists of fixed points of local uniformization group
of the same type. If $s_i=s|_{B_i}: B_i \mapsto E_i$ is a
generalized Fredholm bundle of index $r_i < r$, 
where $E_i$ is the subbundle of $E$ which consists of fixed
points of the local uniformization group,
then in the above proof, one can show
that $e([s:B\mapsto E])$ is in $H_*(B, \ZZ)$.
In the case of Gromov-Witten invariants for rational curves
(cf. section 2, 3), if the target
manifold is semi-positive, then the above assumptions hold. This explains
why the Gromov-Witten invariants in [RT1], [RT2] are integer-valued.

If the above assumptions do not hold, 
in order to get integral Euler classes,
one has to use special properties of certain generalized Fredholm 
orbifold bundles which arise from concrete applications.
\end{rem}

We believe that Theorem 1.1, 1.2 can be 
generalized to other singular varieties.
However, the resulting Euler class may lie in
intersection homology groups.

We end this section with a remark. Let us denote
formally by $\infty_E$ the rank of $E$ and by $\infty _B$
the dimension of $B$. Then $r=\infty_B - \infty _E$.
The Euler class we defined is Poincare dual to
the Chern class $c_{\infty _E}(E)$ (formally) of $E$.
A natural question is whether or not one can construct reasonable cycles,
which are Poincare dual to the lower degree
Chern classes $c_{\infty _E -i}(E)$ of $E$, at least under certain
assumptions on $[s: B \to E]$. The compactness is the main
problem.

\section{Gromov-Witten invariants}
\label{sec:2}

Let $(V, \omega, J)$ be a compact symplectic manifold, where
$\omega$ is a symplectic form and $J$ is a compatible almost
complex structure, i.e., 
\begin{displaymath}
  J^2 = - i d \; , \qquad \omega (J u, J v) = \omega (u, v) \; ,
  \quad \forall u, v \in T V \; .
\end{displaymath}
Then $g (u, v) = \omega (u, J v)$ defines a Riemannian metric on
$V$. Without loss of generality, we may assume that $(V, \omega, J)$ is 
$C^\infty$-smooth.

As usual, if $2g+k \ge 3$, 
we denote by $\M_{g,k}$ the moduli space of
Riemann surfaces of genus $g$ and with $k$ marked points.  Each
point of $\M_{g,k}$ can be presented as
\begin{displaymath}
  \left( \Sigma ; x_1, \ldots, x_k \right)
\end{displaymath}
where $\Sigma$ is a Riemann surface of genus $g$, $x_1, \ldots, x_k
\in \Sigma$ are distinct.  We identify $(\Sigma ; x_1, \ldots,
x_k)$ with $(\Sigma' ; x'_1, \ldots, x'_k)$, if there
is a biholomorphism $f : \Sigma \to \Sigma'$ carrying $x_i$ to
$x_i '$.  Therefore, $\M_{0, 3}$ consists of one point.

Let $\overline{\M}_{g, k}$ be the Deligne-Mumford compactification
of $\M_{g,k}$.  Then $\overline{\M}_{g, k}$ consists of all genus
$g$ stable curves with $k$ marked points.  It is well-known that
$\overline{\M}_{g, k}$ is a K\"{a}hler orbifold (cf. [Mu]).

In this section, we will apply Theorem 1.2 to construct
the GW-invariant
\begin{displaymath}
 \Psi^V_{(A, g, k)} : H^* (V, \QQ)^k \times
      H^* (\overline{\M}_{g, k}, \QQ) \to \QQ \; .
\end{displaymath}

First let us introduce the notion of stable maps.

\begin{defn}
  A stable $C^\ell$-map $(\ell \geq 0)$ with
$k$ marked points is a
tuple $(f, \Sigma; x_1, \ldots, x_k)$ satisfying: 

\renewcommand{\theenumi}{\arabic{enumi})}
\begin{enumerate}
\item \label{item:1}
  $\Sigma = \bigcup^m_{i = 1}
  \sum_i$ is a connected curve with normal crossings and $x_1, \ldots,
  x_k$ are distinct smooth points in $\Sigma$;
\item \label{item:2}
  $f$ is continuous, and each restriction $f |_{\Sigma_i}$
lifts to a $C^\ell$-smooth map from the normalization $\tilde \Sigma_i$ into
$V$;
\item \label{item:3}
  If the homology class of $f |_{\Sigma_i}$ is zero in $H_2(V, \QQ )$
  and $\Sigma _i$ is a smooth rational curve, then $\Sigma_i$ contains at
  least three of $x_i, \ldots, x_k$ and those points in $\Sing
  (\Sigma)$, the latter denotes the singular set of $\Sigma$.
\end{enumerate}
\end{defn}

This definition is inspired by the holomorphic stable maps
in [KM]. We should also note that $2g+k$ may be less than $3$ in the above
definition.

Given $(f, \Sigma; x_1, \ldots, x_k)$ as above, let
$\Aut (\Sigma; x_1, \ldots, x_k)$ be the automorphism group of
$(\Sigma; x_1, \ldots, x_k)$.  Note that if $2g + k \ge 3$, 
$(\Sigma; x_1, \ldots, x_k) \in
\overline{M}_{g, k}$ if and only if 
$\Aut (\Sigma; x_1, \ldots, x_k)$ is
finite. Let ${\rm Aut}(f,\Sigma;x_1, \cdots, x_k)$
be the group consisting of all $\sigma $ in $\Aut (\Sigma; x_1, \ldots, x_k)$
such that $f\cdot \sigma = f$. Clearly, if $(f,\Sigma;x_1, \cdots, x_k)$
is stable, then ${\rm Aut}(f,\Sigma;x_1, \cdots, x_k)$ is finite.

We say that two stable maps $(f, \Sigma; x_1, \ldots,
x_k)$ and $(f', \Sigma'; x'_1, \ldots, x'_k)$ are equivalent if
there is a biholomorphism $\sigma : \Sigma \mapsto \Sigma'$
such that $\sigma (x_i) = x'_i$ $ (1 \leq i \leq k)$ and $f' = f
\circ \sigma$. We will denote by $[f, \Sigma; x_1, \ldots,
x_k]$, usually abbreviated as $[\C]$, 
the equivalence class of stable maps equivalent to $(f,
\Sigma; x_1, \ldots, x_k)$. Note that in case $\Sigma = \Sigma '$,
$\sigma$ is in $\Aut (\Sigma; x_1, \ldots, x_k)$.

The genus of a stable map $(f, \Sigma; x_1, \ldots, x_k)$
is defined to the genus of $\Sigma$. 

Let $\overline {\F}_A^\ell (V, g, k) $ ($\ell \ge 0$) be the space of
equivalence classes $[f, \Sigma; x_1, \ldots, x_k]$
of $C^\ell$ stable maps $ (f, \Sigma; x_1, \ldots, x_k)$
of genus $g$ and with total homology class $A$, which is
represented by the image $f(\Sigma)$ in $V$. 
Clearly, $\overline {\F}_A^\ell (V, g, k)$ is contained 
in $\overline {\F}_A^{\ell'} (V, g, k)$,
if $\ell > \ell'$. 

We will also denote $\overline {\F}_A^0 (V, g, k)$ by
$\overline {\F}_A (V, g, k)$.

For any sequence of stable $C^\ell$-maps 
$\{(f_i, \Sigma_i; x_{i1}, \cdots, x_{ik})\}$, we say that
$(f_i, \Sigma_i; x_{i1}, \cdots, x_{ik})$ converges
to $(f_\infty, \Sigma_\infty; x_{\infty 1}, \cdots, x_{\infty k})$ 
in $C^\ell$-topology, if there are (1) $(\Sigma_i; \{x_{ij}\})$ converges
to $(\Sigma_\infty; \{x_{\infty j}\})$ as marked curves;
(2) $f_i$ converges to $f_\infty$ in $C^0$-topology on $\Sigma_\infty$;
(3) $f_i$ converges to $f_\infty$ in $C^\ell$-topology 
on any compact subset outside the singular set of $\Sigma_\infty$.
Let $\C_i$ be any sequence of equivalence classes of $C^\ell$-stable
maps. We say that $[\C_i]$ converges to $[\C_\infty]$, if 
there are $\C_i=(f_i, \Sigma_i; x_{i1}, \cdots, x_{ik})$ representing
$[\C_i]$ and converging to a representative 
$\C_\infty = (f_\infty, \Sigma_\infty; x_{\infty 1}, \cdots, x_{\infty k})$ 
of $[\C_\infty]$ in $C^\ell$-topology. 

The topology of $\overline {\F}_A^\ell (V, g,k)$ is given by the sequencial
convergence in the above sense. One can easly show that 
the homotopy class of $\overline {\F}_A^\ell (V, g, k)$
is independent of $\ell$.

We define $\F _A(V,g,k)$ to be
the set of all equivalence classes of stable maps with
smooth domain. Put $\F_A^\ell(V,g,k) = \F_A(V,g,k)\cap
\overline{\F}_A^\ell(V,g,k)$.

\begin{rem}
$\F_A^\ell(V,g,k)$ is basically a family of spaces of maps from
Riemann surfaces into $V$. Its topology has been extensively studied in
the literature of algebraic topology. Here, one can regard
$\overline \F_A(V,g,k)$ as a partial compactification of
$\F_A(V,g,k)$. This partial compactification seems to have more
structures than the original space does. The authors do not know
much study on it in the literature. We believe that it 
deserves more attention.
\end{rem}

If $2g + k \ge 3$, one can define
a natural map $\pi_{g,k}$ from $\overline {\F}_A (V, g, k) $
onto $\overline {\M}_{g,k}$ as follows:          

\begin{displaymath}
  \pi_{g,k} ( f, \Sigma; x_1, \ldots, x_k) = 
{\rm Red} (\Sigma; x_1, \ldots, x_k)
\end{displaymath}
where ${\rm Red} (\Sigma; x_1, \ldots, x_k)$ is the stable reduction 
of $(\Sigma; x_1, \ldots, x_k)$, which is obtained by contracting all 
its non-stable irreducible components.
Then, we have $\F _A (V,g,k) = \pi_{g,k}^{-1} (\M _{g,k})$, 
moreover, we can describe
$\F _A^\ell (V,g,k)$ locally as follows: given any $[f,\Sigma;x_1,\cdots, x_k]$
in $\pi_{g,k}^{-1}(\M_{g,k})$. Then the automorphism group 
$\Gamma = {\rm Aut}(\Sigma; x_1,\cdots, x_k)$
is finite. We denote by $\Gamma_0$ its subgroup consisting of automorphisms
preserving $f$. Let
$W_0$ be a small neighborhood of $(\Sigma;x_1,\cdots, x_k)$ in
$\M _{g,k}$, and $p_{W_0}: \tilde W_0 \mapsto W_0$ 
be the local uniformization.
Note that $\Gamma$ acts on $\tilde W_0$ and $W_0= \tilde W_0 / \Gamma $. 
One can show that $[f,\Sigma;x_1,\cdots,x_k]$ has a neighborhood
of the form $\tilde W_0 \times U/\Gamma _0$,
where $U$ is some open neighborhood
of $0$ in the space $C^\ell (\Sigma, f^*TV)$ 
of $f^*TV$-valued $C^\ell$-smooth
functions. Note that $\Gamma_0$ acts on $C^\ell(\Sigma, f^*TV)$ naturally.
Therefore, $\F_A^\ell (V,g,k)$ is a Banach orbifold.

Without much more difficulty, one can also show that 
$\F_A^\ell (V,g,k)$ is a Banach orbifold, even if $2g+k < 3$.
However, it seems to be much harder to prove 
that $\overline {\F}_A^\ell (V, g,k)$ is smooth.
Fortunately, we can avoid it by exploring its weakly smooth
structure. 

Next we define a generalized bundle $E$ over $\overline {\F}_A^1(V, g,k)$.

In the following, we will often denote by $\C$
a stable map $(f, \Sigma; x_1, \ldots, x_k)$,
$f_{\C}$ the map $f$ and $\Sigma_\C$ the connected curve $\Sigma$.

We define $\wedge ^{0,1}_{\C}$ as follows: if
$\Sigma _{\C}$ is smooth, then $\wedge ^{0,1}_{\C}$ consists
of all continuous sections $\nu$ in ${\rm Hom }(T\Sigma _{\C}, f_{\C}^*TV)$
with $\nu \cdot j_\C = - J \cdot \nu$, where
$j_\C$ denotes the complex structure on $\Sigma _\C$.
In other words, $\wedge ^{0,1}_{\C}$ consists of all 
$f_{\C}^*TV$-valued (0,1)-forms $\nu$ over $\Sigma _{\C}$.
In general, $\wedge^{0,1}_\C$ consists of all 
$f_{\C}^*TV$-valued (0,1)-form $\nu$ over 
the normalization of $\Sigma_\C$, more precisely,
if $\Sigma_\C$ has nodes $q_1, \cdots, q_s$, then 
$\wedge ^{0,1}_{\C}$ consists of all 
$f_{\C}^*TV$-valued (0,1)-form $\nu$ over
${\rm Reg}(\Sigma _{\C})$ of $\Sigma _{\C}$ satisfying: 
for each $i$, if $D_1$ and $D_2$ are the two local
components of $\Sigma_\C$ near $q_i$,
then $\nu |_{D_1}$, $\nu |_{D_2}$ can be extended continuously
across $q_i$.
 
Let $\C = (f, \Sigma; \{x_i\})$ and 
$\C'=(f', \Sigma'; \{x_i'\})$ be two equivalent
stable maps, and $\sigma $ be the biholomorphism
from $\Sigma'$ to $\Sigma$ such that $\sigma (x_i') = x_i$
and $f' = f\cdot \sigma$. For convenience,
we sometimes denote $\C'$ by $\sigma^*(\C)$. One can show 
$$\wedge ^{0,1}_{\C'} = \sigma ^* \left (\wedge^{0,1}_{\C}\right ).$$
It follows that $\wedge ^{0,1}_{\C}$ descends to a
space $E_{[\C]}$ over the equivalence class of $\C$.
We put $E= \bigcup _{[\C]} E_{[\C]}$ and equip it 
with the continuous topology. Then $E$ is
a topological fibration over $\overline {\F}_A^1(V,g,k)$. 

For simplicity, we will also use $E$ to denote the restriction
of $E$ to $\overline {\F}_A^\ell (V, g,k)$ for any $\ell > 1$.

There is a natural section $\Phi([\C]): \overline {\F}_A^1(V,g,k) \mapsto E$,
i.e., the Cauchy-Riemann equation, defined as follows:
for any $C^1$-smooth equivalence class 
$[\C]\in \overline {\F}_A^1(V,g,k) $, we
define $\Phi ([\C])$ to be represented by
$$df_{\C} + J \cdot df_{\C} \cdot j_\C \in E_\C,$$
where $j_\C$ denotes the conformal structure of $\Sigma _{\C}$. 
Sometimes, by abusing the notations, we
simply write 
$$\Phi(\C) = df_{\C} + J \cdot df_{\C} \cdot j_\C.$$
Then we have

\begin{proposition}
For any $\ell \ge 2$, 
the section $\Phi : \overline {\F}_A^\ell (V, g, k) \mapsto E$ gives
rise to a generalized
Fredholm orbifold bundle with the natural orientation
and of index $2c_1(V)(A)+2k + (2n-6)(1-g)$.
\end{proposition}

We will postpone the proof of proposition 2.2 to
section 3.

Let $\omega'$ be another symplectic form on $V$ and
$J'$ be one of its compatible almost complex structure, 
Recall that $(\omega', J')$ is deformation equivalent to $(\omega, J)$,
if there is a smooth family of symplectic forms $\omega _s$ and compatible
almost complex structures $J_s$ ($0\le s\le 1$) such that
$(\omega _0, J_0)=(\omega , J)$ and 
$(\omega _1, J_1)=(\omega' , J')$.

\begin{proposition}

Let $\Phi': \overline {\F}_A^\ell (V, g, k) \mapsto E$
be the admissible section induced by the Cauchy-Riemann 
equation of $J'$, where $(\omega', J')$ is given as above. 
Assume that $(\omega', J')$ is deformation equivalent to $(\omega, J)$.
Then $\Phi'$ is homotopic to $\Phi$ as generalized Fredholm orbifold bundles.

\end{proposition}

The proof of this proposition is identical to that of Proposition 2.2.
 
Using the last two propositions, we can construct
symplectic invariants, particularly, the GW-invariants.

In the following, if $2g + k < 3$, for convenience,
we denote by 
$\overline {\M}_{g,k}$ the topological space of one point.

Notice that for any $\ell > 0$,
$\overline {\F}_A^\ell (V, g,k)$ is
homotopically equivalent to $\overline {\F}_A(V, g, k)$.
Then we can deduce from Theorem 1.2 

\begin{theorem}
Let $(V, \omega, J)$ be a compact symplectic manifold with
compatible almost complex structure. Then for each $g, k$ 
and $A \in H_2 (V, \ZZ)$, there is
a symplectically invariant homomorphism
$$\rho^V_{A,g,k} : H^*(\overline {\M}_{g,k}, \QQ) \mapsto
H_*({\overline {\F}}_A (V, g, k), \QQ),$$ 
satisfying: 
for any $\alpha $, $\beta$ in $H^*(\overline {\M}_{g,k}, \QQ)$,
$$\rho ^V_{A,g,k} (\alpha \cup \beta ) = \rho ^V_{A,g,k} 
(\alpha) / \pi_{g,k}^*\beta,$$
where $\pi_{g,k}:{\overline {\F}}_A (V, g, k)
\mapsto \overline {\M}_{g,k}$ is defined as above.
We usually write
$\rho^V_{A,g,k} (1)$ as $e_A(V, g,k)$, which 
is a symplectically invariant class in
$$H_{2c_1(V)(A) + 2k + (2n-6) (1-g)}({\overline 
{\F}}_A (V, g, k), \QQ).$$

Furthermore, if $A=0$, then for any $\beta $ in $\overline \M_{g,k}$, we have
that $\rho ^V_{A,g,k}(\beta)$ takes values in 
$\tau _*(H_*(\overline \M_{g,k}
\times V, \QQ))$,
where $\tau : \overline \M_{g,k} \times V \mapsto \overline F_A(V, g,k)$
is the natural embedding of constant maps. 

\end{theorem}
\vskip 0.1in
\noindent
{\bf Proof:}
By Proposition 2.2, $\Phi: \overline {\F}_A (V, g, k)\mapsto E$
is a genaralized Fredholm orbifold bundle of index $r$,
where $r=2c_1(V)(A) + 2k + (2n-6) (1-g)$.
By Proposition 2.3, its homotopy class 
is independent of choices of $(\omega, J)$.
It follows from Theorem 1.2
that there is an Euler class 
$e([\Phi: \overline {\F}_A (V, g, k)\mapsto E])$
in $H_{r}({\overline 
{\F}}_A (V, g, k), \QQ)$. Then $\rho ^V_{A,g,k}$ is obtained by
taking slant product of this Euler class by cohomological classes
in $H^*(\overline {\M}_{g,k}, \QQ)$. All the properties can be easily 
checked.
\vskip 0.1in
\begin{rem}
We conjecture that the invariant $\rho^V_{g,k}$ is integer-valued, i.e.,
for any $\alpha $ in 
$H^r(\overline {\M}_{g,k}, \ZZ)$, $\rho ^V_{g,k} (\alpha)$
is in $H_{2c_1(V)(A) + 2k + 2n (1-g)-r}({\overline 
{\F}}_A (V, g, k), \ZZ)$.
\end{rem}

In order to define the GW-invariants, we
observe that there is an evaluation map
$$
\begin{array}{rl}
ev : \overline {\F}_A (V, g, k) &\mapsto V^k, \\
ev(f, \Sigma; x_1, \ldots , x_k) &= (f(x_1), \ldots, f(x_k)),\\
\end{array}
$$
then we can define the $GW$-invariants
\begin{displaymath}
 \Psi^V_{(A, g, k)} : H^* (V, \QQ)^k \times
      H^* (\overline{M}_{g, k}, \QQ) \to \QQ \; ,
\end{displaymath}
namely, for any $\alpha_1, \ldots, \alpha_k \in H^* (V, \QQ)$,
$\beta \in H^* (\overline{M}_{g, k}, \QQ)$,
\begin{displaymath}
 \Psi^V_{(A, g, k)} (\beta; \alpha_1, \ldots, \alpha_k) =
  \ev^* 
  \left( \pi_1^* \alpha_1 \wedge \cdots \wedge \pi^*_k \alpha_k
  \right) (\rho ^V_{A,g,k}(\beta ))
\end{displaymath}

\begin{theorem}

  The $GW$-invariants $\Psi^V_{(A, g, k)}$ are symplectic
  invariants satisfying the composition law and the basic
  properties {\rm (cf. [KM], [RT2], [T], [Wi1])}.
\end{theorem}

\begin{rem}
  It is believed that $\Psi^V_{(A, g, k)}$ is also
 integer-valued. In fact, it is true for semi-positive symplectic
manifolds (cf. [RT1], [RT2]).
 \end{rem}

\begin{exam}
  Let $(V, \omega, J)$ be a symplectic manifold as above
and $\omega$ be an integer class.
  Then for any holomorphic curve $C \subset V$,
  \begin{displaymath}
    \int_C \omega \geq 1
  \end{displaymath}
  We say that a pseudo-holomorphic map $f : S^2 \to V$ is a line
  if $\int_{S^2} f^*\omega  = 1$.  Let $A$ be the homology class of
  lines, then the moduli space of lines is compact modulo
  automorphisms of $S^2$.

  On the other hand, $G = \Aut (S^2)$ acts naturally on 
  $\Map_A (S^2, V)$ and the bundle $\wedge^{0, 1} (T V)$ over 
  $\Map _A(S^2, V)$.  Therefore, we have a Fredholm bundle $E$ over
  $B = \Map_A (S^2, V)_0 / G$, where $\Map_A (S^2, V)_0$ denotes the space
  of maps which are generically immersive.  The Cauchy-Riemann
  equation descends to a section of $E \to B$.
One can show that 
  \begin{displaymath}
 \Psi^V_{(A, 0, 3)} (\alpha_1, \alpha_2, \alpha_3) =
  \Bigl(\ev \bigl( \pi^{-1} (e (B, E))  \bigr)
  \cap (\alpha^*_1 \times \alpha^*_2 \times \alpha^*_3) \Bigr)
  \hbox{~in~} V^3
  \end{displaymath}
  where $\pi: \Map_A(S^2, V) \mapsto B$ is the natural
projection, and 
$\alpha^*_1 \times \alpha^*_2 \times \alpha^*_3$ are
  Poincare duals of $\alpha_1, \alpha_2, \alpha_3$.

  Note that $e (B, E) \in H_r (B, \ZZ)$ for $r = 2 (c_1 (V) \cdot
  A + n - 3)$.
\end{exam}

We end up this section with two basic decomposition properties
of the symplectic invariant $\rho ^V_{A,g, k}$.

Let $\sigma : \overline \M_{g_1, k_1 +1} \times \overline \M_{g_2, k_2 +1}
\mapsto \overline \M_{g, k}$, where $g = g_1 + g_2$ and $k = k_1 + k_2$
with $2g_1+k_1\ge 2$, $2g_2+k_2\ge 2$,
be the map by glueing the $k_1+1$-th marked point of the first factor
to the first marked point of the second factor. We denote by
$PD(\sigma)$ the Poincare dual of ${\rm Im}(\sigma)$.
The composition law expresses $\rho ^V_{A,g,k} (PD(\sigma))$
in terms of $\rho ^V_{A_1,g_1, k_1+1}$ and $\rho ^V_{A_2,g_2, k_2 +1}$
with $A=A_1+A_2$.

Given any decomposition $A= A_1 + A_2$, there is 
an natural map 
$$\begin{array}{rl}
&p: \overline \F_{A_1}(V, g_1, k_1 +1)\times 
\overline \F_{A_2}(V, g_2, k_2 +1)
\mapsto V\times V,\\
&p([h_1, \Sigma_1; x_{1}, \cdots, x_{k_1+1}], [h_2, \Sigma_2;
y_1,\cdots, y_{k+2+1}]) = (h_1(x_{k_1+1}), h_2(y_1)).\\
\end{array}
$$
Let $\Delta$ be the diagonal in $V\times V$. Then
there is an obvious map $\pi$ from
$p^{-1}(\Delta)$ onto $\pi_{g,k}^{-1}({\rm Im}(\sigma))$ by
identifying $x_{k_1+1}$ with $y_1$.

Clearly, $\rho ^V_{A,g,k} (PD(\sigma))$ can be regarded as a
class in $H_*(\pi_{g,k}^{-1}({\rm Im}(\sigma)), \QQ)$.

On the other hand, if $\{u_i\}$ is any basis of $H^*(V, \ZZ)$ and 
$\{u^*_i\}$ is its dual basis, then we have a homology class
$$\sum _i \rho ^V_{A_1, g_1, k_1+1}/
ev^* \pi_{k_1 +1}^* u_i \otimes \rho ^V_{g_2, k_2+1}/ ev^*\pi_1^* u_i^*$$
in $H_*(p^{-1}(\Delta), \QQ)$.

The first composition law for $\rho^V_{A,g,k}$ is given by
the equation:
$$\rho ^V_{A,g,k} (PD(\sigma))= \pi_*\left (\sum _{A=A_1+A_2}
\sum _i \rho ^V_{A_1, g_1, k_1+1}/
ev^* \pi_{k_1 +1}^* u_i \otimes 
\rho ^V_{A_2,g_2, k_2+1}/ ev^*\pi_1^* u_i^* \right ).
$$

The second composition law for $\rho ^V_{A,g,k}$
arises from the map $\theta : \overline \M_{g-1, k+2}
\mapsto \overline \M_{g, k}$, which is obtained by glueing the
last two marked points, in a similar way. 

As above, we define
$$\begin{array}{rl}
&p: \overline \F_{A}(V, g-1, k +2)
\mapsto V\times V,\\
&p([h, \Sigma; x_{1}, \cdots, x_{k+2}]) = (h(x_{k+1}), h(x_{k+2})).\\
\end{array}$$
We also have the resolution
$\pi: p^{-1}(\Delta) \mapsto \pi_{g,k}^{-1}({\rm Im}(\theta))$. 
Then we have
$$\rho ^V_{A,g, k}(PD(\theta)) = \pi_*\left (\sum _i\rho ^V_{A,g-1, k+2}/ \pi_{k+1}^*u_i
\wedge \pi^*_{k+2} u_i^*\right ).$$

\section{The proof of Proposition 2.2 and 2.3}
\label{sec:3}

In the section, we prove Proposition 2.2 in details.
The same arguments can be applied to proving Proposition 2.3.
We will omit its proof except a few comments at the end of this section.

Fix any $\ell \ge 2$. Let
$(f, \Sigma; x_1,\ldots, x_k)$ be a stable $C^\ell$-map representing
a point in $\overline {\F}_A^\ell (V, g, k) $.
Since the structure of $\F_A^\ell (V,g,k)$ is clear (cf. section 2), we
may assume that $[f, \Sigma; x_1,\ldots, x_k]$
is in $\overline {\F}_A^\ell (V, g, k) \backslash 
\overline {\F}_A^\ell (V, g, k) $.

The components of $\Sigma$ can be grouped into two parts:
the principal part and the bubbling part.
The principal part consists of those components 
of genus bigger than zero and those rational components, which contain 
at least three of $x_1, \cdots, x_k$ and the points in 
${\rm Sing} (\Sigma )$.
Other non-stable rational components consist in the bubbling
part. 

By adding one or two
marked points to each bubbling component,
we obtain a stable curve $(\Sigma; x_1, \cdots, x_k, z_1,\ldots, z_l)$
in $\overline{\M}_{g, k +l}$, where $z_1, \ldots, z_l$ are
added points.

Let $W$ be a small neighborhood of $(\Sigma; \{x_i\},\{ z_j\})$ 
in $\overline{\M}_{g, k + l}$, and $\tilde W$ be
the uniformization of $W$. Then $W=\tilde W/ \Gamma$, where
$\Gamma = {\rm Aut}(\Sigma; \{x_i\}, \{z_j\})$.

If $2g + k \ge 3$, we can express $\tilde W$ as follows:
by contracting the bubbling part,
we obtain the stable reduction $(\Sigma'; y_1, \ldots, y_k)$ 
of $(\Sigma; x_1, \ldots , x_k)$.
Let $W_0$ be a small neighborhood of $(\Sigma'; y_1, \ldots,
y_k)$ in $\overline{\M}_{g,k}$. 
Let $\tilde W_0$ be the uniformization of $W_0$.
Then $W_0=\tilde W_0/ \Gamma$ and
$\tilde W = W\times _{W_0}\tilde W_0$. 
In particular, $\tilde W = W$ is smooth whenever $W_0$ is smooth.

Let $\tilde \U $ be the universal family of curves over $\tilde W$. 
Clearly, $\tilde \U$ is smooth. We fix a metric $g$ on $\tilde \U$. 
For any two maps $h_1$, $h_2$ from fibers of $\U$ over $\tilde W$,
we define the distance 
 \begin{eqnarray*}
 d_{\tilde W}(h _1, h_2) = & \sup _{x\in {\rm Dom (h_1)}}
\sup _{d_g (y, x) = d_g (x,  {\rm Dom (h_2)})}
d_V (h_1 (x), h_2(y))\\
+~&\sup _{y\in {\rm Dom (h_2)}}
\sup _{d_g (x, y) = d_g (y,  {\rm Dom (h_1)})}
d_V (h_1 (x), h_2(y)),\\
  \end{eqnarray*}
where 
$d_g(\cdot, \cdot )$, $d_V(\cdot , \cdot )$ are
distance functions of $g$ and $V$.

Let $\Sigma_j$ be any non-stable
component of $(\Sigma; x_1, \ldots, x_k)$, then
by the definition, the homology class of $f(\Sigma_j)$
is nontrivial. It follows that there is at least one regular value
for $f|_{\Sigma_j}$.
Therefore, we may choose
$z_1, \ldots , z_{l}$, such that for each
$i$, $f^{-1} (f(z_i))$ consists of
finitely many immersive points. 

Choose local hypersurfaces $H_1, \ldots, H_{l}$, such that
each $H_i$ intersects ${\rm Im}(f)$ transversally at $f(z_i)$.

Fix a small $\delta > 0$, we define
$$
\begin{array}{rl}
\Map _\delta (W) = &\{ (\tilde f, \tilde
\Sigma; \{\tilde x_i\},\{ \tilde z_j\})~ |~
(\tilde \Sigma; \{\tilde x_i\}, \{\tilde z_j\}) \in \tilde W,   
d_{\tilde W} (\tilde f, f) < \delta,\\
&~~~\tilde f ~{\rm is }~C^0 ~{\rm on~}\tilde \Sigma {\rm ~and}~ C^\ell
\end{array}
$$
We will equip it with the topology: any sequence $(h_a, \Sigma_a; \{x_{ia}\},
\{ z_{ja}\} )$ converges to $(h_\infty, \Sigma_{\infty}; \{x_{i\infty}\},
\{ z_{j\infty}\})$, if $(\Sigma_a; \{x_{ia}\}, \{ z_{ja}\})$ converges to
the stable curve $(\Sigma_{\infty}; \{x_{i\infty}\}, \{ z_{j \infty}\})$ 
in $\tilde W$, and $h_a$ converges to $h_\infty$ in $C^0$-topology 
everywhere and $C^\ell$-topology outside $\Sing (\Sigma_\infty)$.

We denote by $\Sing (\tilde \U)$ the union of singularities of
the fibers of $\tilde \U$ over $\tilde W$.
Let $K$ be any compact subset in $\tilde \U \backslash \Sing(\tilde \U)$
of the form: there is a diffeomorphism $\psi _K: (K\cap \Sigma)
\times \tilde W\mapsto K$ such that $\psi_K ((K\cap \Sigma )\times \{t\})$
lies in the fiber of $\tilde \U$ over $t$. Then 
we define
$$\begin{array}{rl}
&\Map _\delta (W, K)\\
=&\{ (\tilde f, \tilde
\Sigma; \{\tilde x_i\}, \{\tilde z_j\} )\in \Map_\delta (W)~|~
||(\tilde f \cdot \psi_K  - f)|_{K\cap \Sigma \times \{0\}} ||_{C^\ell } 
< \delta \}.\\
\end{array}
$$
Clearly, each $Map _\delta (W, K)$ is open in $Map _\delta (W)$.

By forgetting added marked points, each point in
$\Map _\delta (W)$ gives rise to a stable map
$\C$ and consequently, an equivalence class $[\C]$ in 
$\overline {\F}_A^\ell (V, g, k)$. 
Let $p_W$ be such a projection map into $\overline {\F}_A^\ell (V, g, k)$.
We denote by $\Map _\delta(W_0, K)$ the image of $\Map _\delta (W,K)$ in 
$\overline {\F}_A^\ell (V, g, k)$ under the projection $p_W$.

Let ${\rm Aut}(\C)$ be the automorphism group 
of the stable map $\C$. It is a
subgroup of $\Gamma$, so it is finite and acts on $\tilde \U$.
Let us denote by $m(\C)$ its order.

>From now on, $K$ always denotes a compact set in $\tilde \U \backslash
\Sing (\tilde \U)$ containing an open
neighborhood of $\bigcup _j f^{-1}(f(z_j))$.
Moreover, we may assume that $K$ is invariant under the action of ${\rm Aut}(\C)$.
\begin{lemma}
If $\delta > 0$ is sufficiently small, then
the map $p_W|_{\Map _\delta (W,K)}$ is finite-to-one of the order $m(\C)$, 
and $\Map _\delta (W_0, K)$ is an open
neighborhood of $\C$ in $\overline {\F}_A^\ell (V, g, k)$.
  
Furthermore, there is a canonical action of ${\rm Aut}(\C)$
on $\Map_\delta (W,K)$ with the quotient $\Map_\delta (W_0, K)$.
In particuar, if $\C$ has trivial automorphism group, then 
$p_W|_{\Map _\delta (W,K)}$ is actually one-to-one.
\end{lemma}

\vskip 0.1in
\noindent
{\bf Proof:} Suppose that  
$(h', \Sigma'; \{x_i'\}, \{z_j'\})$ and $(h'', \Sigma''; \{x_i''\}, 
\{z_j''\})$, which are close to $\C$, 
have the same image under the projection $p_W$.
Then there is a biholomorphism $\sigma: \Sigma ' \mapsto \Sigma ''$,
such that $h' = h''\cdot \sigma$ and $\sigma (x'_i) = x_i''$.
Since $h', h''$ are close to $f$,
$\sigma $ has to be close to an automorphism of $\C$.
Since $h'(z_j'), h''(z_j'') \in H_j$ for $1\le j \le l$,
we have that $\sigma (z'_j)$ and $z_j''$ are close
to $f^{-1}(f(z_j))$. Since $f$ is transverse to
$H_j$ for each $j$, $p_W$ is finite-to-one
of order no more than $ m(\C)$.

Let us construct the action of ${\rm Aut}(\C)$ on
$\Map_\delta (W,K)$ with $\Map _\delta (W_0, K)$ as its quotient.
In fact, let $\tau \in {\rm Aut}(\C)$ and 
$\C'=(h', \Sigma'; \{x_i'\}, \{z_j'\})$ in $\Map _\delta (W,K)$.
If $\C'$ is very close to $\C$, then there is a unique
sequence $\{z_{\tau j}\}$ in $\tau (\Sigma')$ 
such that $h'(\tau ^{-1}(z_{\tau j}) )\in
H_j$ and $z_{\tau j}$ is very close to $z_j$. We put 
$$\tau_*(\C') = (h'\cdot \tau^{-1}, \tau (\Sigma'); \{\tau(x_i')\}, \{
z_{\tau j}\}),$$
then $\tau_*(\C) \in \Map _\delta (W,K)$. 
Clearly, if $\tau'$ is another one in ${\rm Aut}(\C)$,
we have that 
$$(\tau\cdot \tau')_*(\C') = \tau _* (\tau'_*(\C')).$$
It follows that
there is an natural action
of ${\rm Aut}(\C)$ on $\Map _\delta (W,K)$.
Clearly, the quotient is $\Map _\delta (W_0,K)$. It also follows
that $p_W$ is of order $m(\C)$.

If $(h, \Sigma'; x_1, \ldots, x_k)$ is a stable map 
very close to $\C$, then $h$ is immersive near $z_j$
and there are unique $z_j'$ in $\Sigma '$ near $z_j$, such that
$h(z_j') \in H_j$. It follows that $(h, \Sigma'; \{x_i'\}, \{z_j'\})$
is in $\Map _\delta (W, K)$, so $\Map_\delta (W_0,K)$ is
a neighborhood of $[\C]$ in $\F^\ell _A(V,g,k)$.
The lemma is proved.
\vskip 0.1in

Recall that a $TV$-valued, $(0,1)$-form over
the universal family $\tilde \U$ of curves
is a continuous section $\nu$ in $Hom (\pi_1^*T\tilde \U, \pi^*_2TV)$
satisfying: $\nu \cdot j_{\tilde \U} = - J\cdot \nu$, where
$j_{\tilde \U}$ denotes the complex structure on $\tilde \U$.
We denote by $\Gamma ^{0,1}_\ell (\tilde \U, TV)$ the space of
such $(0,1)$-forms, which are $C^\ell$ smooth and vanish
near ${\rm Sing}(\tilde \U)$.

Note that $E|_{\Map_\delta (W_0, K)} \mapsto \Map _\delta (W_0, K)$
lifts to a topological bundle, denoted by $E|_{\Map_\delta (W, K)}$,
or simply $E$ if no possible confusions, 
over $\Map _\delta (W,K)$.
In order to prove Proposition 2.2, we need to show 
that each $E|_{\Map_\delta (W, K)}$
is a generalized Fredholm bundle over $\Map_\delta (W, K)$.

Let $\Phi$ be defined by the Cauchy-Riemann equation in section 2.
It lifts to a section, still denoted by $\Phi$, of $E$ over
$ \Map _\delta (W)$, explicitely, 
$$\Phi (\tilde f, \tilde \Sigma; \{
\tilde x_i\}, \{\tilde z_j\}) = 
d\tilde f + J \cdot d\tilde f \cdot j_{\tilde \Sigma}.
$$
Let $L_{\tilde f}$ be the linearization of $\Phi$ 
at $\tilde f$. Then,
for any vector field $u$ over $\tilde f (\tilde \Sigma)$,
$$L_{\tilde f}(u) = du + J(\tilde f) \cdot du \cdot j_{\tilde \sigma} + \nabla _u J
\cdot d\tilde f \cdot j_{\tilde \Sigma}.$$

We denote by $r$ the distance function
to the singular set $\Sing (\tilde \U)$ with respect to $g$.

For any smooth section $u \in \Gamma ^0(\tilde \Sigma , \tilde f^*TV)$, 
we define the norm 
$$||u||_{1,p} = 
\left (
\int _{\tilde \Sigma} (|u|^p + |\nabla u|^p ) d\mu \right )^{\frac{1}{p}}
 + \left ( \int _{\tilde \Sigma}
r^{-\frac{2(p-2)}{p}} |\nabla u|^2 d\mu \right )^{\frac{1}{2}},  
$$
where $p \ge 2$ and $\Gamma ^0(\tilde \Sigma , \tilde f^*TV)$ 
is the space of continuous
sections of $\tilde f^*TV$ over $\tilde \Sigma $, and all norms, covariant
derivatives are taken with respect to the metric $g|_{\tilde \Sigma}$.
If $\tilde \Sigma$ has more than one components, then $u$ consists of
continuous sections of components which have the same value at each node.

Then we define
$$
\begin{array}{rl}
L^{1,p}(\tilde \Sigma, \tilde f^*TV)=
\{ u \in \Gamma ^0(\tilde \Sigma , \tilde f^*TV) ~|~
||u||_{1,p} < \infty \}.
\end{array}
$$

\begin{lemma}

For any $p \ge 2$, there is a uniform constant $c(p)$ 
such that for any fiber $\tilde \Sigma$ 
of $\tilde \U$ over $\tilde W$, and any $u$ in
$L^{1,p}(\tilde \Sigma, \tilde f^*TV)$, we have 
$$||u||_{C^0} \le c(p) ||u||_{1,p}.$$

\end{lemma}

\vskip 0.1in
\noindent
{\bf Proof:} We observe that any small geodesic ball of $\tilde \Sigma$
is uniformly equivalent to an euclidean ball or the union of
two euclidean annuli of the same size.
Then the lemma follows from the standard Sobolev Embedding Theorem.
\vskip 0.1in

It follows that
$L^{1,p} (\tilde \Sigma, \tilde f^*TV)$ is complete for $p > 2$. 

On the other hand, for any 
$v \in {\rm Hom} (T\tilde \Sigma, \tilde f^*TV)$, we define
$$||v||_{p} = \left (\int _{\tilde \Sigma} |v|^p d\mu 
\right )^{\frac{1}{p}}
+ \left ( \int _{\tilde \Sigma} r ^{-\frac{2(p-2)}{p}}
|v|^2 d\mu \right )^{\frac{1}{2}},
$$
where all norms and derivatives are taken with respect to 
$g|_{\tilde \Sigma}$, too. Then we put
$$ 
L^p (\wedge ^{0,1}(\tilde f^*TV))
=\{ v \in {\rm Hom}(T\tilde \Sigma, \tilde f^*TV) ~|~J\cdot v = - 
v \cdot j_{\tilde \Sigma},~
||v||_p < \infty\}.
$$

For any $(\tilde f, \tilde \Sigma; \{
\tilde x_i\}, \{\tilde z_j\})$ in $\Map _\delta (W)$,
$L_{\tilde f}$ maps the space $L^{1,p}(\tilde \Sigma, \tilde f^*TV)$
into $L^p(\wedge ^{0,1} (\tilde f^*TV))$.
Let $L_{\tilde f}^* $ be its adjoint operator
with respect to the $L^2$-inner product on
$L^2(\wedge ^{0,1} (\tilde f^*TV))$, more explicitly,
for any $\tilde f^*TV$-valued (0,1)-form $v$,
$$L^*_{\tilde f} v = - e_1(v_1)- e_2(v_2) + B_{\tilde f}( v),$$
where $\{e_1, e_2\}$ is any orthonormal basis of $\tilde \Sigma$
with $j_{\tilde \Sigma} (e_1) = e_2$, $v_i = v (e_i)$ ($i= 1,2$)
and $B_{\tilde f}( v)$ is an operator of order $0$, defined
by
$$2 g_V(u, B_{\tilde f} (v)) = g_V((\nabla _u J )e_2(\tilde f), v_1)
- g_V((\nabla _u J)e_1(\tilde f), v_2)$$
for any $u \in L^{1,2} (\tilde \Sigma, \tilde f^*TV)$.
We denote by ${\rm Coker}(L_{\tilde f})$
the space of all $v$ in $L^2(\wedge ^{0,1} (\tilde f^*TV))$ such
that $L^*_{\tilde f} (v) = 0$. Then by the standard elliptic theory,
it is a finite dimensional subspace in $L^p(\wedge ^{0,1} (\tilde f^*TV))$
for any $p$.

\begin{lemma}
For any $v \in L^p(\wedge ^{0,1} (\tilde f^*TV))$ ($p \ge 2$), there are
$v_0 \in {\rm Coker}(L_{\tilde f})$ and 
$u\in L^{1,p}(\tilde \Sigma, \tilde f^* TV)$, such that
$L_{\tilde f} u = v - v_0$.
\end{lemma}

\vskip 0.1in
\noindent
{\bf Proof:} By the definition, one can find
$u\in L^{1,2}(\tilde \Sigma, \tilde f^* TV)$ and $v_0 \in
{\rm Coker}(L_{\tilde f})$, such that
$L_{\tilde f} u = v- v_0$. Then the lemma follows from the standard elliptic theory.
\vskip 0.1in

Let $\C$ be the fixed holomorphic, stable $C^\ell$-map,
in particular, $\Phi (\C) =0$. 

For any $v \in \Gamma ^{0,1}_{\ell -1} (\U, TV)$,
we define its restriction $v|_{\tilde \C}$
to a stable map $\tilde \C$ as follows:
let $\tilde \C = (\tilde f, \tilde \Sigma; \{
\tilde x_i\}, \{\tilde z_j\})$, then for any
$x\in \tilde \Sigma$, we define
$$v|_{\tilde \C} (x) = v (x, \tilde f(x)).$$

Let $S$ be any finitely dimensional subspace in
$\Gamma ^{0,1}_{\ell -1} (\tilde \U, TV)$ ($\ell \ge 2$).
We define 
$$S|_\C = \{ v|_\C ~|~ v \in S\}.$$
Then we can define
$E_S$ over $\Map _\delta (W, K)$ as follows:
for any $\tilde \C$ in $\Map _\delta (W, K)$,
$E_S |_{\tilde \C} = S|_{\tilde \C}$.

Assume that $\dim S = \dim S|_\C$. Then if $K$ in 
$\tilde \U \backslash {\rm Sing}(\tilde \U)$ is 
sufficiently large and $\delta$ is sufficiently small,
$E_S$ is a bundle of rank $\dim S$ over $\Map _\delta (W, K )$.

The following is the main technical result of this section.

\begin{proposition}
Let $S$ be as above. Suppose that
its restriction $S|_\C$ to
$\C$ is transverse to $L_{f_\C}$, i.e.,
if $v_1, \cdots, v_s$ span $S$, then  
$v_1|_\C ,\cdots, v_s|_\C$
and ${\rm Im }(L_{f_\C})$ generate
$L^{1,p}(\wedge ^{0,1}f^*TV)$. 
Then by shrinking $W$ if necessary, if $\delta $ is sufficiently small
and $K$ is sufficiently large,
$\Phi ^{-1}(E_S)$ is a smooth submanifold,
which contains $\C$, in $\Map_\delta (W, K)$ and 
of dimension $2 c_1(V)(A) + 2(n-3)(1-g)+ 2k + \dim S$.
Moreover, $E_S \mapsto \Phi^{-1}(E_S)$ is a smooth bundle.
\end{proposition}

\begin{rem}

Suppose that $W$ and $S$ are invariant under the natural action of ${\rm Aut}(\C)$.
Clearly, there is an induced action of ${\rm Aut}(\C)$
on both $\Phi^{-1}(E_S)$ and the total space
of the bundle $E_S$ over $\Phi^{-1}(E_S)$. 

\end{rem}

Now let us prove Proposition 3.4. The tool is the Implicit 
Function Theorem.

Let $\C$ be the stable
map $(f, \Sigma; x_1, \cdots, x_k)$ in ${\rm Map }(W)$ 
as given in Proposition 3.4.
We denote by $q_1, \cdots, q_s$ those nodes in $\Sigma$.
Recall that $z_1, \cdots, z_l$ be the added points such that $f(z_i)
\in H_i$.

Fix an $\nu$ in $S$ such that its restriction to $\Sigma_\C$
is $0$. In fact, for proving Proposition 3.4, we suffice to take $\nu =0$.

First we want to construct a family of approximated
$(J, \nu_t)$-maps $\tilde f_t$ parametrized by 
$t \in \tilde W$, where
$\nu_t=\nu |_{\Sigma_t}$. 
Note that a $(J, \nu_t)$-map 
is a smooth $\tilde f: \Sigma _t\mapsto V$ satisfying
the inhomogeneous Cauchy-Riemann equation:
$$\Phi(\tilde f)(y) = \nu (y, \tilde f(y)), ~~~~y \in \Sigma_t.$$

For any $q_i$ ($1\le i \le s$), by shrinking 
$\tilde W$ if necessary, we may choose
coordinates $w_{i1}, w_{i2}$, as well as $t$ in $\tilde W$, near $\C$,
such that the fiber 
$$(\Sigma _t; x_1(t), \cdots, x_k(t),
z_1(t), \cdots, z_l(t))$$ 
of $\tilde \U $ over $t$ is locally given by
the equation
$$ w_{i1} w_{i2} = \epsilon _i (t),~~~|w_{i1} | < 1, ~|w_{i2} | < 1, $$
where $\epsilon _i$ is a $C^\infty$-smooth function of $t$.
For any $y$ in $\Sigma _t$, if $|w_{i1} (y)| > L \sqrt {|\epsilon (t)|}$
or $|w_{i2} (y)| > L \sqrt {|\epsilon (t)|}$ for all $i$, where $L$ is a 
large number,
then there is a unique $\pi_t(y) $ in $\Sigma = \Sigma _0$ such that
$d_g(y, \pi_t(y)) = d_g(y, \Sigma)$. Note that if $y$ is not in 
the coordinate chart given by $w_{i1}, w_{i2}$, then we simply set
$w_{i1}(y) = w_{i2}(y) = \infty$.

The following lemma follows from straightforward computations.

\begin{lemma}
 
  For any $k > 0$, there is a uniform constant $a_k$ such that for
any $y$ in $\Sigma _t$ with either $|w_{i1} (y)| > L \sqrt {|\epsilon (t)|}$
or $|w_{i2} (y)| > L \sqrt {|\epsilon (t)|}$ for all $i$,
$$|\nabla ^k (\pi_t - id |_{\Sigma_t})| (y) \le 
a_k \min _i\{ |t|, \frac{|\epsilon_i(t)|}{ d_g(y, q_i)^{k+1}}\},$$
where $\nabla$ denotes the covariant derivative with respect to
$g$, and both $\pi_t$, $id$ are regarded as maps from $\Sigma_t$ into
$\tilde \U$.

\end{lemma}

Let us introduce a complex structure $\tilde J$
on $\tilde \U \times V$ as follows:
for any $u_1 \in T\tilde \U \subset T(\tilde \U\times V)$, 
$$\tilde J (u_1) = j_\U (u_1) + \nu (j_\U(u_1));$$
For any $u_2$ in $TV \subset T(\tilde \U\times V)$, we put
$\tilde J (u_2) = J(u_2)$.

Define $F: \Sigma \mapsto \tilde \U\times V$ by assigning $y$ in $\Sigma$
to $(y, f(y))$. We call $F$ the graph map of $f$.
One can show that $F$ is $\tilde J$-holomorphic. 
In fact, for any given $\tilde f: \Sigma _t \mapsto V$, it is
a $(J, \nu_t)$-map if and only if its graph map 
is $\tilde J$-holomorphic (cf. [Gr]).

Put $p_i = F(q_i)$. Without loss of generality, 
we may assume that 
$$F(\{ w_{i1} w_{i2} = 0 ||w_{i1}| < 1, |w_{i2}| < 1\} )$$
is contained in a coordinate chart $(u_1, \cdots, u_{2N})$ 
of $\tilde \U \times V$ near $p_i$. 
We may further assume that
$$
\begin{array}{rl}
\tilde J ( \frac{\partial~}{\partial u _i} ) &= 
\frac{\partial~~}{\partial u_{N+i}} + {\cal O} (|u| ),~\\
\tilde J ( \frac{\partial~~}{\partial u_{N+i}} ) &= - 
\frac{\partial~}{\partial u_i} + {\cal O} (|u| ),\\
\end{array}
$$
where $i = 1,2, \ldots, N$ and $ | u | = \sqrt{\sum_{i=1} ^{2n} |u_i |^2 }$.
The curve $F(\Sigma)$ has two components near $p_i$, 
which intersect transversally there. Then by
changing $u_1, \cdots , u_{2N}$ appropriately, we may assume that in
complex coordinates $u_1 + \sqrt{-1} u_{N+1}$, ..., $u_N + \sqrt{-1} u_{2N}$,
$$ 
F(w_{i1}, w_{i2}) = (w_{i1}, w_{i2} , 0 ,0,\cdots, 0 ) 
+ {\cal O} ( | w_{1i} | ^2 + |w_{i2}|^2 ).
$$
Using this same formula, one can easily extend 
$F$ to a neighborhood of $q_i$ in $\U$.

\begin{lemma}

Let $\pi_2 :\tilde  \U \times V \mapsto V$ be the natural projection. Then
there is a uniform constant $a$ such that for any $y$ in $\Sigma_t$
with $\frac{1}{ 2} \le |w_{i1}(y)| \le 1$ or $\frac{1}{2}
\le |w_{i2}(y)| \le 1$,
$$|\pi_2(F(y)) - f(\pi_t(y))|_{C^2}\le a |\epsilon _i (t)|.$$

\end{lemma}

This lemma can be easily proved by straightforward computations.

Let $\eta: \RR^1 \mapsto \RR^1 $ be a cut-off function
satisfying: $\eta (x) = 0$ for $|x |\le 1$, $\eta (x) =1$ for
$|x| > 2$, and $|\eta ^{(k)}(x) | \le 2^k$. 

We define
$\tilde f_t(y)$, where $y \in \Sigma _t$, as follows: if either 
$|w_{i1} (y)| > 1$
or $|w_{i2} (y)| > 1$ for all $i$,
put $\tilde f_t (y) = f(\pi_t(y))$; If for some $i$,
$|w_{i1} (y)| < \frac{1}{2} $ and 
$|w_{i2} (y)| < \frac{1}{2}$, then
we define $\tilde f_t(y) = \pi_2(F(y))$;
If $\frac{1}{2} \le |w_{i1} (y)| \le 1$ or
$\frac{1}{2} \le |w_{i2} (y)| \le 1$, we define
$\tilde f_t(y)$ to be 
$${\rm exp}_{f(q_i)} \left ( 
\eta (2 d_g(y, q_i)) 
{\rm exp}_{f(q_i)}^{-1}
f(\pi_t(y)) + (1 - \eta (2 d_g(y, q_i)) )
{\rm exp}_{f(q_i)}^{-1} \pi_2(F(y)) \right ).
$$

Since $f$ is continuous at each $q_i$, $\tilde f_t$ is
continuous.

\begin{lemma}

There is a uniform constant $a_f$ such that for any $0 \le k \le \ell$
and $y \in \Sigma _t $,
$$
\begin{array}{rl}
|\nabla ^k \tilde f_t|(y) &\le a_f \min \{|t|,
   \frac{|\epsilon(t)|}{d_g(y, q_i)^{k+1}}\},\\
|\nabla ^{k-1}(\Phi (\tilde f_t)- \nu(\cdot, \tilde f_t(\cdot)))|(y)
   &\le a_f \min _i\{ |t|, \frac{|\epsilon_i(t)|}{d_g(y, q_i)^{k}}\}.\\
\end{array}
$$
In particular, $|\nabla \tilde f_t|$ is uniformly bounded.

\end{lemma}

\vskip 0.1in
\noindent
{\bf Proof:} By Lemma 3.5 and 3.6, we suffice to prove those
estimates near a given node, say $q_i$. Assume that
$|w_{i1}|(y), |w_{i2}|(y) < \frac{1}{2}$.
Then $\tilde f_t(y) = \pi_2(F(y))$.
Let us prove the second estimate.
The proof for the first is identical. We omit it. 

We may assume that $|w_{i1}(y)| \ge |w_{i2}(y)|$. 
Let $J_0$ be the standard complex structure in the coordinate chart 
$\{u_1, \cdots u_{2N}\}$, i.e.,
$$J_0(\frac{\partial }{\partial u_i})=\frac{\partial }{\partial u_{N+i}},
J_0(\frac{\partial }{\partial u_{N+i} })= - \frac{\partial }{\partial u_i},$$
where $i = 1, \cdots, N$. Note that $F$ is holomorphic with respect to 
$J_0$. Then we can deduce
$$
\begin{array}{rl}
 &\Phi (\tilde f_t)(y) - \nu (y, \tilde f_t(y))\\ 
= &d\pi_2 \cdot \left (dF + 
\tilde J \cdot dF\cdot j_{\Sigma_t} \right ) 
(w_{i1}(y), w_{i2}(y), 0, \cdots,0)\\
=& d\pi_2 \cdot (\tilde J - J_0)
\cdot dF \cdot j_{\Sigma_t }(w_{i1}(y), w_{i2}(y), 0, \cdots,0)\\
\le & c |\nabla F| |\tilde J - J_0|(w_{i1}(y), w_{i2}(y), 0, \cdots,0),\\
\end{array}
$$
where $c$ is some unfiorm constant.
It follows that 
$$|\Phi (\tilde f_t ) (y)- \nu(y, \tilde f_t(y)| \le \frac{c |\epsilon
(t)|} {d_g(y, q_i)}.$$
Similarly, one can
deduce other cases of the second estimate from the above identity.
\vskip 0.1in

For any $t$ small, we denote by $g_t$ the induced metric
on $\Sigma_t$ by $g$. Note that $r$ is the distance function
from $\Sing (\tilde \U)$ with respect to $g$.

For any smooth section $u \in \Gamma ^0(\Sigma _t, \tilde f_t^*TV)$, we
recall 
$$||u||_{1,p} = 
\left (
\int _{\Sigma_t} (|u|^p + |\nabla u|^p ) d\mu_t \right )^{\frac{1}{p}}
 + \left ( \int _{\Sigma_t}
r^{-\frac{2(p-2)}{p}} |\nabla u|^2 d \mu_t \right )^{\frac{1}{2}},  
$$
and
$$
\begin{array}{rl}
L^{1,p}(\Sigma _t, \tilde f_t^*TV)=
\{ u \in \Gamma ^0(\Sigma _t, \tilde f_t^*TV) ~|~
||u||_{1,p} < \infty \},
\end{array}
$$
where $p \ge 2$ and $\Gamma ^0(\Sigma _t, \tilde f_t^*TV)$ 
is the space of continuous
sections of $\tilde f_t ^*TV$ over $\Sigma _t$.
If $\Sigma_t$ has more than one components, then $u$ consists of
sections which are continuous over each of its components 
and have the same value at each node.

We put 
$$L^{1,p}=
\{ (u,t)~|~u \in L^{1,p}(\Sigma _t, \tilde f_t^*TV)\}.
$$
It is a topological bundle over $\U$.

On the other hand, for any $v \in {\rm Hom} (T\Sigma_t,\tilde f_t^*TV)$, 
we have
$$||v||_{p} = \left (\int _{\Sigma_t} |v|^p d\mu_t \right )^{\frac{1}{p}}
+ \left ( \int _{\Sigma_t} r ^{-\frac{2(p-2)}{p}}
|v|^2 d\mu_t \right )^{\frac{1}{2}},
$$
and
$$
\begin{array}{rl} 
L^p (\wedge ^{0,1}(\tilde f_t^*TV))
=\{ v \in {\rm Hom}(T\Sigma_t, \tilde f^*TV) ~|~J\cdot v = - 
v \cdot j_{\Sigma_t},~
||v||_p < \infty\}.
\end{array}
$$
As above, we put $L^p (\wedge ^{0,1}(TV))$ to be the union of
all $L^p (\wedge ^{0,1}(\tilde f_t^*TV))$ with $t \in \tilde W$. It is
another topological bundle over $\U$.

Furthermore, if $C^\ell_0 (\tilde \U, TV)$
denotes the space of all $C^\ell$-smooth sections,
which vanish near ${\rm Sing}(\tilde \U)$, of $\pi^*_2TV$ over
$\U\times V$,
then there is an embedding of $C^\ell_0(\tilde \U, TV)$ into
$L^{1,p}$, where $\pi_2: \tilde \U\times V \mapsto V$ is the natural projection. 

Similarly, there is an embedding of $\Gamma^{0,1}_{\ell -1}(\tilde \U, TV)$
into $L^p(\wedge ^{0,1}TV)$.

Note that both $C^\ell _0(\tilde \U, TV)$ and 
$\Gamma^{0,1}_{\ell -1}(\tilde \U, TV)$ are bundles over $\tilde \U$.

By straightforward computations, we
can deduce from Lemma 3.7.

\begin{lemma}
For any $p > 2$, we have 
$$||\Phi (\tilde f_t )-\nu(\cdot, \tilde f_t(\cdot))||_p 
\le c |t|^{\frac{1}{2}},$$
where $c$ is a uniform constant.
\end{lemma}

Next we define a map from $L^{1,p}$ into $L^{p} (\wedge ^{0,1}(TV))$
as follows: for any $(u, t)$ in $L^{1,p}$, 
$$\Psi (u, t) = \Phi ( exp _{\tilde f_t} u ),$$
where $exp_{\tilde f_t} u$ denotes the function which takes value
$exp _{\tilde f_t(x)} u(x)$
at $x$. 

Clearly, this map $\Psi$ is well-defined, and
maps $C^\ell_0 (\tilde \U, TV)$ into $\Gamma^{0,1}_{\ell -1}(\tilde \U, TV)$.

\begin{rem}
Here we have used the fact that $(x,t) \mapsto (x,t, \tilde f_t(x))$
defines a smooth map from $\tilde \U$ into $\tilde \U\times V$.
\end{rem}

Now let us study the linearization
$L_t = D_u\Psi$ of $\Psi$ at $(0, t)$:
for any $u$, we have 
$$L_t(u) =  du + J(\tilde f_t) \cdot du \cdot j_{\Sigma_t} 
+ \nabla _u J (\tilde f_t) \cdot d \tilde f_t \cdot 
j_{\Sigma_t}.
$$
First we want to establish uniform elliptic estimates for 
$L_t$.

\begin{lemma}

There is a uniform constant $c$ such that for any $(u, t)$
in $L^{1,p}$, we have
$$||u||_{1,p} \le c ( ||L_t(u) ||_p + ||u||_{1,2}).$$

\end{lemma}

\vskip 0.1in
\noindent
{\bf Proof:}
We may assume that $p > 2$, otherwise, the lemma is trivially true.
Without loss of generality, we may further assume that
$r\le  \frac{1}{2}$ if both $w_{i1}$ and $ w_{i2}$
are less than $ \frac{1}{2}$ for some $i$.

Let $\eta$ be a cut-off function satisfying:
$\eta(x) = 0$ for $|x| \le \frac{1}{4}$, $\eta(x) = 1$ for $|x| > \frac{1}{2}$,
and $|\eta'| \le 2$.

Put $\tilde u = \eta (r) u$. It vanishes whenever
$|w_{i1}|, |w_{i2}| \le \frac{1}{2}$. Moreover, we have
$$L_t(\tilde u) = \eta (r) L_t(u) + \eta'(r) \left (
u \,dr + (J(\tilde f_t) u )\,dr\cdot j_{\Sigma_t} \right ).$$

Since $\Sigma _t$ has uniformly bounded geometry in the region where
$r \ge \frac{1}{4}$, we can apply the standard $L^p$-estimate
for $1^{st}$-order elliptic operators and obtain
$$||\tilde u||_{1,p} \le c (||L_t(\tilde u)||_p + ||\tilde u||_{1,2}).$$
Together with the previous identity, we deduce
$$||\tilde u||_{1,p} \le c (||L_t(u)||_p + ||u||_p + ||u||_{1,2}).$$
Note that $c$ always denotes a uniform constant, which may depend on $p$.
By the Sobolev inequality in dimension two ( $\dim \Sigma_t =2$),
we have $||u||_p \le c ||u||_{1,2}$, hence,
$$||\tilde u||_{1,p} \le c (||L_t(u)||_p + ||u||_{1,2}).$$

Therefore, we suffice to show that
for each $i$,
$$
\begin{array}{rl}
\left ( \int _{|w_{i1}|, |w_{i2}| \le  \frac{1}{2}}
|\nabla u |^p d\mu_t \right )^{\frac{1}{p}}
&\le c (||L_t(u)||_p + ||u||_{1,2}),\\
\left ( \int _{|w_{i1}|, |w_{i2}| \le \frac{1}{2}}
(|w_{i1}|^2 + |w_{i2}|^2)^{-\frac{p-2}{p}}|\nabla u|^2d\mu_t
\right )^{\frac{1}{2}}
&\le c (||L_t(u)||_p + ||u||_{1,2}).\\
\end{array}
$$
Let us first prove the second inequality. 
Without loss of generality, we assume that $\epsilon _i(t) \not= 0$.
Write
$w_{i1} = \rho e^{\sqrt{-1} \theta}$, then
$w_{i2} = \frac{|\epsilon_i(t)|}{\rho} e^{\sqrt{-1} (\theta + \theta_0)}$,
where $\epsilon _i(t) = |\epsilon_i(t)| e^{\sqrt{-1} \theta_0}$. 
Hence, $|w_{i1}|^2 + |w_{i2}|^2 = 
\rho ^2 + \frac{|\epsilon_i(t)|^2}{\rho^2}$.
Moreover, $|w_{i1}|, |w_{i2}| \le 1$ whenever 
$|\epsilon_i(t)| \le \rho \le 1$.

Put $u_i$ to be zero if either
$\rho > 1$ or $\rho < |\epsilon_i(t)|$, and $ (1-\eta (\frac{r}{2})) u$ otherwise.
In terms of $\rho$ and $\theta$, we have the following expression:
$$L_t(u_i) (\frac{\partial }{\partial \rho})
= \frac{\partial u_i }{\partial \rho} + \frac{1}{\rho}
J(\tilde f_t)\left (
\frac{\partial u_i}{\partial \theta}\right ) +\frac{1}{\rho}
 (\nabla _{u_i} J )\frac{\partial
\tilde f_t}{\partial \theta}.$$

It follows that
$$
\begin{array}{rl}
&\int _{|w_{i1}|, |w_{i2}| \le 1}
(|w_{i1}|^2 + |w_{i2}|^2)^{-\frac{p-2}{p}}
|\frac{\partial u_i }{\partial \rho} + \frac{1}{\rho}
J(\tilde f_t)\left (
\frac{\partial u_i}{\partial \theta}\right )|^2d\mu_t\\
\le & c \left ( ||L_t(u_i)||_p^2 + 
\int _{|w_{i1}|, |w_{i2}| \le 1}
(|w_{i1}|^2 + |w_{i2}|^2)^{-\frac{p-2}{p}} |u_i|^2 |\nabla \tilde f_t|^2
d\mu_t \right )\\
\le & c \left (||L_t(u)||_p^2 + ||u||_{1,2}^2 + 
\int _{|w_{i1}|, |w_{i2}| \le 1}
(|w_{i1}|^2 + |w_{i2}|^2)^{-\frac{p-2}{p}} |u|^2
d\mu_t \right )\\
\end{array}
$$
Notice that the integral
$$\int _{|w_{i1}|, |w_{i2}| \le 1}
(|w_{i1}|^2 + |w_{i2}|^2)^{-\frac{p-2}{p-1}}
d\mu_t$$
is bounded by a constant depending only on $p$.
However, by the Sobolev Embedding Theorem,
we have
$$\left (\int _{|w_{i1}|, |w_{i2}| \le 1}
|u_i|^{2p} d\mu_t \right )^{\frac{1}{p}}
\le c(p) ||u||_{1,2}^2.$$ 
It follows
$$
\begin{array}{rl}
&\int _{|w_{i1}|, |w_{i2}| \le 1}
(|w_{i1}|^2 + |w_{i2}|^2)^{-\frac{p-2}{p}}
|\frac{\partial u_i }{\partial \rho} + \frac{1}{\rho}
J(\tilde f_t)\left (
\frac{\partial u_i}{\partial \theta}\right )|^2d\mu_t\\
\le & c ( ||L_t(u)||_p^2 + ||u||_{1,2}^2).\\
\end{array}
$$

We have
$$
\begin{array}{rl}
&\int _{|w_{i1}|, |w_{i2}| \le 1}
(|w_{i1}|^2 + |w_{i2}|^2)^{-\frac{p-2}{p}}
\left ( |\frac{\partial u_i }{\partial \rho}|^2 + \frac{1}{\rho^2}
|\frac{\partial u_i}{\partial \theta}|^2\right ) d\mu_t\\
=&\int _{|w_{i1}|, |w_{i2}| \le 1}
(|w_{i1}|^2 + |w_{i2}|^2)^{-\frac{p-2}{p}}
\left (|\frac{\partial u_i }{\partial \rho} + \frac{1}{\rho}
J(\tilde f_t)\left (
\frac{\partial u_i}{\partial \theta}\right )|^2\right .\\
&~~~~~~\left .
- 2 \langle \frac{\partial u_i }{\partial \rho},
\frac{1}{\rho}
J(\tilde f_t)\left (
\frac{\partial u_i}{\partial \theta}\right )\rangle 
\right )  d\mu_t\\
\end{array}
$$
Using integration by parts, we derive 
$$
\begin{array}{rl}
&\left | \int _{|w_{i1}|, |w_{i2}| \le 1}
(|w_{i1}|^2 + |w_{i2}|^2)^{-\frac{p-2}{p}}
 \langle \frac{\partial u_i }{\partial \rho},
\frac{1}{\rho} J_0 \left (
\frac{\partial u_i}{\partial \theta}\right )\rangle 
d \mu_t \right | \\
\le & \frac{p-2}{p}\int _{|w_{i1}|, |w_{i2}| \le 1}
(|w_{i1}|^2 + |w_{i2}|^2)^{-\frac{p-2}{p}}
| \langle u_i- a(\rho), \frac{1}{\rho^2} J_0 \left (
\frac{\partial u_i}{\partial \theta}\right )\rangle |
d \mu_t,\\
\end{array}
$$
where $J_0 = J(q_i)$ and $a(\rho)$ is any function on $\rho$. 
Using the Poincare inequality on the unit
circle and choosing $a(\rho)$ appropriately, we can show that
the last integral is no bigger than
$$\int _{|w_{i1}|, |w_{i2}| \le 1}
(|w_{i1}|^2 + |w_{i2}|^2)^{-\frac{p-2}{p}}
\frac{1}{\rho^2} |\frac{\partial u_i}{\partial \theta}|^2
d \mu_t,$$
which is the same as the integral
$$
\begin{array}{rl}
&\frac{1}{2} \int _{|w_{i1}|, |w_{i2}| \le 1}
(|w_{i1}|^2 + |w_{i2}|^2)^{-\frac{p-2}{p}}\\
&~~\left ( |\frac{\partial u_i}{\partial \rho}- L_t(u_i)- \frac{1}{\rho}
(\nabla _{u_i} J) \frac{\partial \tilde f}{\partial \theta}|^2+
\frac{1}{\rho^2} |\frac{\partial u_i}{\partial \theta}|^2 \right )
d \mu_t\\
\le &  c\, (||u_i|_{1,2}^2 + ||L_tu_i||^2_p)\\
&~+ \left ( \frac{1}{2} + \frac{1}{4p} \right )
\int _{|w_{i1}|, |w_{i2}| \le 1}
(|w_{i1}|^2 + |w_{i2}|^2)^{-\frac{p-2}{p}}
\left ( |\frac{\partial u_i}{\partial \rho}) |^2+
\frac{1}{\rho^2} |\frac{\partial u_i}{\partial \theta}|^2 \right )
d \mu_t.\\
\end{array}
$$
On the other hand, by the above arguments,
one can also show that
$$
\begin{array}{rl}
&\left | \int _{|w_{i1}|, |w_{i2}| \le 1}
(|w_{i1}|^2 + |w_{i2}|^2)^{-\frac{p-2}{p}}
 \langle \frac{\partial u_i }{\partial \rho},
\frac{1}{\rho} (J-J_0) \left (
\frac{\partial u_i}{\partial \theta}\right )\rangle 
d \mu_t \right |\\
\le & c ||u_i||_{1, 2} + \frac{1}{2p} 
\int _{|w_{i1}|, |w_{i2}| \le 1}
(|w_{i1}|^2 + |w_{i2}|^2)^{-\frac{p-2}{p}}
\left ( |\frac{\partial u_i}{\partial \rho} |^2+
\frac{1}{\rho^2} |\frac{\partial u_i}{\partial \theta}|^2 \right )
d \mu_t .\\
\end{array}
$$
Combining all the above inequalities, we can deduce
the second inequality we wanted.

To obtain the first from the second, we decompose
the region $\{ |\epsilon_i(t)| \le \rho \le 1$ into
subannuli $\{ \delta _j \le \rho \le \delta _{j-1}\}$, where
$j=1, \cdots, m$, $\delta_0 =1$, $\delta_m = |\epsilon_i(t)|$ and
$1\le \frac{\delta _{j-1}}{\delta_j} < 2$.

On each subannulus $\{ \delta _j \le \rho \le \delta _{j-1}\}$,
the scaled metric $\delta^{-2}_j g_t$ has bounded geometry, we
can apply the standard $L^p$-estimate and obtain

$$
\begin{array}{rl}
& \int _{ \delta _j \le \rho \le \delta _{j-1}}
|\nabla u_i|^p d\mu_t   ~~~~~~~~~~~~~~~~\\
\le c& \left ( \int_{ \delta _j \le \rho \le \delta _{j-1}}
|L_tu_i|^p d\mu_t
+ \left (\delta_j^{-\frac{2p-4}{p}} \int_{ \delta _j \le \rho \le \delta _{j-1}}
|\nabla u_i|^2 d\mu_t \right )^{\frac{p}{2}}\right ).\\
\end{array}
$$

Clearly, the first inequality we wanted 
follows by suming up these over $j$.
The lemma is proved. 
\vskip 0.1in

\begin{lemma}
Let $S$ be as in Proposition 3.4
and $t$ sufficiently 
small. Then for any $p > 2$ and $v$ in 
$L^p (\wedge ^{0,1} \tilde f_t ^*TV)$, there are 
$u$ in $L^{1,p}(\tilde f_t^* TV)$ and $v_0$ in $S$,
satisfying:
$$\begin{array}{rl}
& L_t u  =v - v_0,\\
&\max \{||u||_{1,p}, ||v_0||_p\} \le c ||v||_p,\\
\end{array}
$$
where $c$ is a uniform constant. 
\end{lemma}

\noindent
{\bf Proof:} First we prove that there are $u$, $v_0$ such that
$L_t u = v -v_0$ for sufficiently small $t$. If not, we can find 
a sequence $\{t_j\}$ with $\lim t_j =0$ and $v_j$ in ${\rm Coker}(L_{t_j})$,
such that each $v_j$ is perpendicular to $S$ with respect to
the $L^2$-metric on $L^2 (\wedge ^{0,1} \tilde f_{t_j} ^*TV)$.

Note that $v_j \in L^p (\wedge ^{0,1} \tilde f_{t_j} ^*TV)$ for any $p$. We normalize
$||v_j||_p = 1$.

By using standard elliptic estimates (cf. [GT]), one can easily show
that $v_j$ converges to some $v_\infty $ in $L^p (\wedge ^{0,1} \tilde f_\C ^*TV)$
outside the singular set of $\Sigma _\C$. Clearly, $v_\infty$ is perpendicular
to $S$ and $L_0^*v_\infty =0$, so by our assumptions, $v_\infty = 0$. It follows
that for any compact subset $K' \subset \U$ with $K' \cap {\rm Sing}(\Sigma _\C)
= \emptyset$, we have
$$\int _{K' \cap \Sigma_{t_j}} |v_j|^p d\mu_{t_j} 
+ \int _{K' \cap \Sigma_{t_j}} r ^{-\frac{2(p-2)}{p}}
|v_j|^2 d\mu_{t_j}  \mapsto 0, ~~{\rm as}~j\mapsto \infty.
$$

Put $t=t_j$ for any fixed $j$.
Let $\epsilon_i (t), w_{i1}, w_{i2}$ be as above, near some
node $q_i$ of $ \Sigma _\C$.
As before, without loss of generality, we assume that $\epsilon _i(t) \not= 0$.
Write $w_{i1} = \rho e^{\sqrt{-1} \theta}$, then
$w_{i2} = \frac{|\epsilon_i(t)|}{\rho} e^{\sqrt{-1} (\theta + \theta_0)}$,
where $\epsilon _i(t) = |\epsilon_i(t)| e^{\sqrt{-1} \theta_0}$. 
Hence, $|w_{i1}|^2 + |w_{i2}|^2 = 
\rho ^2 + \frac{|\epsilon_i(t)|^2}{\rho^2}$.
Moreover, $|w_{i1}|, |w_{i2}| \le 1$ whenever 
$|\epsilon_i(t)| \le \rho \le 1$.

Using $L_{t_j}^* v_j = 0$, we have
$$
\begin{array}{rl}
&\int _{|w_{i1}|, |w_{i2}| \le 1}
(|w_{i1}|^2 + |w_{i2}|^2)^{-\frac{p-2}{p}}
|\frac{\partial v_{j\rho } }{\partial \rho} + \frac{1}{\rho}
\frac{\partial v_{j\theta}}{\partial \theta}|^2d\mu_t\\
\le & c  
\int _{|w_{i1}|, |w_{i2}| \le 1}
(|w_{i1}|^2 + |w_{i2}|^2)^{-\frac{p-2}{p}} |v_j|^2
d\mu_t,\\
\end{array}
$$
where $v_{j\rho } = v_j (\frac{\partial }{\partial \rho })$
and $v_{j\theta } = v_j (\frac{\partial }{\partial \theta})$.
Note that 
$$v_{j\rho } = J(\tilde f_t) v_{j \theta},~~~
v_{j\theta } = - J(\tilde f_t) v_{j \rho}.$$
It follows that
$$
\begin{array}{rl}
&\int _{|w_{i1}|, |w_{i2}| \le 1}
(|w_{i1}|^2 + |w_{i2}|^2)^{-\frac{p-2}{p}}
|\frac{\partial v_{j\rho } }{\partial \rho} - \frac{1}{\rho}
J(\tilde f_t)\left (
\frac{\partial v_{j \rho }}{\partial \theta}\right )|^2d\mu_t\\
\le & c \int _{|w_{i1}|, |w_{i2}| \le 1}
(|w_{i1}|^2 + |w_{i2}|^2)^{-\frac{p-2}{p}} |v_j|^2
d\mu_t.\\
\end{array}
$$

We have
\begin{eqnarray*}
\int _{|w_{i1}|, |w_{i2}| \le 1}
(|w_{i1}|^2 + |w_{i2}|^2)^{-\frac{p-2}{p}}
\left ( |\frac{\partial v_{j\rho} }{\partial \rho}|^2 + \frac{1}{\rho^2}
|\frac{\partial v_{j\rho}}{\partial \theta}|^2\right ) d\mu_t\\
=\int _{|w_{i1}|, |w_{i2}| \le 1}
(|w_{i1}|^2 + |w_{i2}|^2)^{-\frac{p-2}{p}}
\left (|\frac{\partial v_{j\rho} }{\partial \rho} - \frac{1}{\rho}
J(\tilde f_t)\left (
\frac{\partial v_{j\rho}}{\partial \theta}\right )|^2\right .\\
\left .
+ 2 \langle \frac{\partial v_{j\rho} }{\partial \rho},
\frac{1}{\rho}
J(\tilde f_t)\left (
\frac{\partial v_{j\rho}}{\partial \theta}\right )\rangle 
\right )  d\mu_t
\end{eqnarray*}
Using integration by parts, we derive 
$$
\begin{array}{rl}
&\int _{|w_{i1}|, |w_{i2}| \le 1}
(|w_{i1}|^2 + |w_{i2}|^2)^{-\frac{p-2}{p}}
 \langle \frac{\partial v_{j\rho} }{\partial \rho},
\frac{1}{\rho} J_0 \left (
\frac{\partial v_{j\rho}}{\partial \theta}\right )\rangle 
d \mu_t \\
=&\int _{|w_{i1}|=1~{\rm or}~|w_{i2}|=1}
(|w_{i1}|^2 + |w_{i2}|^2)^{-\frac{p-2}{p}}
 \langle v_{j\rho} ,
\frac{1}{\rho} J_0 \left (
\frac{\partial v_{j\rho}}{\partial \theta}\right )\rangle 
d \mu_t\\
+& \frac{p-2}{p}\int _{|w_{i1}|, |w_{i2}| \le 1}
(|w_{i1}|^2 + |w_{i2}|^2)^{-\frac{p-2}{p}}
 \langle v_{j\rho}- a(\rho), \frac{1}{\rho^2} J_0 \left (
\frac{\partial v_{j\rho}}{\partial \theta}\right )\rangle 
d \mu_t,\\
\end{array}
$$
where $J_0 = J(q_i)$ and $a(\rho)$ is any function on $\rho$. 
Using the Poincare inequality on the unit
circle and choosing $a(\rho)$ appropriately, we can show that
the last integral is no bigger than
$$\int _{|w_{i1}|, |w_{i2}| \le 1}
(|w_{i1}|^2 + |w_{i2}|^2)^{-\frac{p-2}{p}}
\frac{1}{\rho^2} |\frac{\partial v_{j\rho}}{\partial \theta}|^2
d \mu_t.$$
On the other hand, 
one may assume that for $j$ sufficiently large,
$$
\begin{array}{rl}
&\int _{|w_{i1}|, |w_{i2}| \le 1}
(|w_{i1}|^2 + |w_{i2}|^2)^{-\frac{p-2}{p}}
 \langle \frac{\partial v_{j\rho} }{\partial \rho},
\frac{1}{\rho} (J-J_0) \left (
\frac{\partial v_{j\rho}}{\partial \theta}\right )\rangle 
d \mu_t \\
\le & \frac{1}{2p} 
\int _{|w_{i1}|, |w_{i2}| \le 1}
(|w_{i1}|^2 + |w_{i2}|^2)^{-\frac{p-2}{p}}
\left ( |\frac{\partial v_{j\rho}}{\partial \rho} |^2+
\frac{1}{\rho^2} |\frac{\partial v_{j\rho}}{\partial \theta}|^2 \right )
d \mu_t .\\
\end{array}
$$
Combining all above estimates, we have
$$\begin{array}{rl}
&\lim _{j \to \infty} \int _{|w_{i1}|, |w_{i2}| \le 1}
(|w_{i1}|^2 + |w_{i2}|^2)^{-\frac{p-2}{p}} |\nabla v_{j}|^2
d \mu_t\\
=~&\lim_{j \to \infty}
\int _{|w_{i1}|, |w_{i2}| \le 1}
(|w_{i1}|^2 + |w_{i2}|^2)^{-\frac{p-2}{p}}
|v_{j}|^2
d \mu_t\\
=~&0.\\
\end{array}
$$
Then one can deduce from this that $\lim_{j\to \infty} ||v_j||_p=0$.
A contradiction! Therefore, we have proved the first part.

Let us prove the estimate by contradiction. Suppose that it is
not true, then there are $u_i$ in $L^{1,p}(\tilde f_{t_i}^*TV)$
and $v_{0i}$ in $S$ satisfying: 

\noindent
(1) $\max\{||u_i||_{1,p},||v_{0i}||_p\}=1$;

\noindent
(2) $u_i$ are perpedicular to ${\rm Ker}(\pi_S\cdot L_{t_i})$,
where $\pi_S$ is the projection onto the orthogonal complement of $S$
in $L^{1,2}(\tilde f_{t_i}^*TV)$;

\noindent
(3) $\lim_{i\to \infty} ||L_{t_i}u_i + v_{0i}||_p =0$.

We may choose $t_i$ such that $\lim _i t_i = t_\infty$ exists.

By (1) and the Sobolev Embedding Theorem, 
we may assume that $u_i$ converges to $u_\infty$
in the $L^{1,2}$-norm. We may further assume that
$v_{0i}$ converges to some $v_{0\infty}$.
Note that $L_{f} u_\infty = v_{0\infty}$.

If $v_{0\infty} \not= 0$, then $u_\infty \not= 0$.
Then $u_\infty \in {\rm Ker}(\pi_S\cdot L_{f})$,
which is impossible. Therefore, we have $v_{0\infty} =0$.
This implies that $\lim ||u_i||_{1,p} = 1$. It follows from Lemma 3.9
that $||u_i||_{1,2}$ is uniformly bounded away from zero.
Then one can show that $u_\infty$ 
is in ${\rm Ker}(\pi_S\cdot L_f)$, a contradiction!
The lemma is proved.
\vskip 0.1in

Let $P$ be a finitely dimensional subspace in $C^\ell_0(\tilde \U, TV)$.
Then for any map $\tilde f: \tilde \Sigma \mapsto V$, 
where $\tilde \Sigma$ is a fiber of $\tilde \U$ over $\tilde W$,
we define 
$u|_{\tilde f}$ by
$$u|_{\tilde f}(x) = u(x, \tilde f(x)), ~~{\rm for ~any~} x\in \tilde \Sigma,$$
and 
$$P_{\tilde f} ~=~\{ u|_{\tilde f} ~|~ u \in P\}.$$
We assume that $\dim P = \dim P_f$ and $q_S({\rm Ker}(\pi_S\cdot L_0))= P_f$,
where $\pi_S$ is defined in the proof of Lemma 3.10
and $q_S: L^{1,2}(\tilde \Sigma, \tilde f^*TV) \mapsto
P_{\tilde f}$ is the projection with respect to
the $L^2$-inner product.

One can easily deduce from the above lemma the following.

\begin{lemma}
Let $P$ and $S$ be as above
and $t$ be sufficiently small. Then for any $p > 2$, $u_0
\in P_{\tilde f_t}$ and $v$ in 
$L^p (\wedge ^{0,1} \tilde f_t ^*TV)$, there are unique 
$u$ in $L^{1,p}(\tilde f_t^* TV)$ and $v_0$ in $S$,
satisfying:
$$\begin{array}{rl}
&q_S(u) =u_0,~~ L_t (u )  = v - v_0,\\
&\max \{||u||_{1,p}, ||v_0||_p\} \le c \max \{||u_0||_{1,p}, ||v||_p\},\\
\end{array}
$$
where $c$ is a uniform constant. 
\end{lemma}

\vskip 0.1in
\noindent
{\bf Proof of Proposition 3.4:} We have the following expansion:

$$\Psi (u, t) = \Psi (0, t) + L_t u + H_t (u), $$
where $H_t(u)$ is the term of higher order satisfying:
$||H_t(u)||_p \le c ||u||_{C^0} ||u||_{1,p}$
for some uniform constant $c$, which depends only on
the derivatives of $J$. 
By the Sobolev Embedding Theorem, it follows
$$||H_t(u)||_p \le c ||u||_{1,p}^2.$$
Also note that $\Psi(0, t) = \Phi (\tilde f_t)$.

Consider the map $\Xi: L^{1,p}\times E_S \mapsto L^p(\wedge ^{0,1}TV)\times E_P$,
defined by
$$\Xi(u,t,v_0) = (\Psi(u,t) + v_0, q_S(u)).$$
Note that $E_P$ is the bundle induced by $P$ over $\tilde W$ with fibers
$P_{\tilde f_t}$.
The linearization of $\Xi$ at $(0,t,0)$ is the map 
$$\begin{array}{rl}
D\Xi : L^{1,p}(\Sigma_t,\tilde f_t^*TV)\times S_{\tilde f_t} 
&\mapsto ~L^p(\wedge ^{0,1}\tilde f_t^*TV)\times P_{\tilde f_t},\\
(u,v_0) ~~~&\mapsto ~~(L_t(u) + v_0 , q_S(u)).\\
\end{array}
$$
By Lemma 3.11, it is an isomorphism 
with uniformly bounded inverse. Therefore, by Lemma 3.8 and
the Implicit Function Theorem, there is an $\epsilon _0> 0$
such that for any $(0,u_0)\in L^p(\wedge ^{0,1}\tilde f_t^*TV)\times P_{\tilde f_t}$
with $||u_0||_{1,p}< \epsilon_0$
and $d_{\tilde W}(t, 0)< \epsilon_0$,
there is a unique $(u,t,v_0)$ satisfying:
$$\begin{array}{rl}
&\Xi(u,t,v_0)= (0, u_0),~~~\\
&\max \{||u||_{1,p}, ||v_0||_p\} \le c ||u_0||_{1,p},\\
\end{array}
$$
where $c$ is some uniform constant. 

It follows that if $W$ is sufficiently small, the subset
$$\{(u,t)\in L^{1,p} | \pi_S\cdot \Psi(u,t)=0, ||u||_{1,p}< \epsilon_0\}$$
is parametrized by $u_0$ in $P$ and $t\in \tilde W$. In particular,
it is a smooth manifold of dimension $\dim S + 2c_1(V)(A)
+2n(1-g)+2k+2l$. Note that by our choice of $P$, we have
$$\dim P = \dim S + 2c_1(V)(A)+ 2n(1-g).$$
We define $Y_{\epsilon_0}(S,W)$ to be
$$\{(u,t)\in L^{1,p} | \pi_S\cdot \Psi(u,t)=0, ||u||_{1,p}< \epsilon_0,
{\rm exp}_{\tilde f_t(z_j)} u(z_j) \in H_j \},$$
where $z_j$ ($1\le j\le l$) are added points
given at the beginning of this section. Then $Y_{\epsilon_0}(S,W)$
is a smooth manifold of dimension
$$\dim S + 2c_1(V)(A) + 2n(1-g)+2k.$$
 
We claim that for $\delta$ sufficiently small and $K$ is sufficiently large,
$\Phi^{-1}(E_S)$ is an open set in $Y_{\epsilon_0}(S,W)$.

Let $(\tilde f, \tilde \Sigma; \{\tilde x_i\}, \{\tilde z_j\})$ be in 
$\Phi^{-1}(E_S)$. We denote by $t$ the corresponding
point $(\tilde \Sigma; \{\tilde x_i\}, \{\tilde z_j\})$ in $\tilde W$.

Using the fact that $d(\tilde f, f) \le \delta$, we can write
$\tilde f(x)= {\rm exp}_{\tilde f_t(x)}u(x)$ for 
some $\tilde f_t^*TV$-valued function $u$.
We suffice to show that $u \in L^{1,p}(\Sigma_t, \tilde f_t^*TV)$
and $||u||_{1,p}< \epsilon_0$.
It follows from the following lemma.

\begin{lemma}
For any $p>2$, there is a uniform constant $c$
such that
$$\int_{\tilde \Sigma} r^{\frac{2(p-2)}{p}} |u|^p d\mu
\le c ||u||_{C^0(K)},$$
where $r$ is the distance function from the set of nodes as we used before.
\end{lemma}

Write $v_0=\Phi(\tilde f)\in S$. Then $||v_0||_{C^0(K)}\le c \delta$
for some uniform constant $c$.
By our choice of $S$, it follows that if $\delta$ is sufficiently small,
$||v_0||_{C^1}<<\epsilon_0$. Then
Lemma 3.12 can be proved by asymptotic analyses near nodes of $\Sigma$
or the arguments in the proof of Lemma 3.9.

Finally, by differentiating $\pi_S\cdot \Phi (\tilde f)=0$ on $t$ and
using Lemma 3.11, one can show that $E_S \mapsto \Phi^{-1}(E_S)$ is
a smooth bundle and $\Phi|_{\Phi^{-1}(E_S)}$ is a smooth section.
This is essentially the smooth dependence of solutions, which are
produced by the Implicit Function Theorem, 
on parameters.
 
Proposition 3.4 is proved.

\vskip 0.1in
\noindent
{\bf Proof of Proposition 2.2:} We first need to construct a covering
of $\Phi^{-1}(0)$ by open subsets, which will be
parametrized by $[\C ]= [f, \Sigma; \{x_i\}] \in \Phi^{-1}(0)$, 
a small number
$\delta > 0$, an neighborhood $W_0$ of the stable reduction
${\rm Red}(\Sigma; \{x_i\})$ of $(\Sigma; \{x_i\})$ in $\overline {\M}_{g,k}$,
a compact subset $K$ in the universal family 
$\tilde \U$ of curves over $\tilde W$.
Here $W$, $\tilde W$ are given as before.
We define 
$$U_\delta ([\C], W_0, K ) = \Map_\delta (W_0, K),$$ 
where $\Map_\delta (W_0, K)$ is given
in Lemma 3.1. 

Each $U_\delta ([\C], W_0, K)$ is of the form
$\Map_\delta (W, K) / \Gamma $, where $\Gamma = {\rm Aut}(\C)$ and 
$\Map_\delta (W, K)$ were given as in Lemma 3.1.
We put 
$$\tilde U_\delta ([\C], W_0, K ) = \Map_\delta (W, K).$$
It is the uniformization of $U_\delta ([\C], W_0, K)$.
Therefore, we have shown that $\overline {\F}^\ell _A (V, g,k)$ is
a topological orbifold.

Let $E$ be the space of $TV$-valued $(0,1)$-forms defined
in section 2. For each $U_\delta ([\C], W_0, K)$, as we have 
already seen, $E$ can be lifted to a topological bundle
$E|_{\tilde U_\delta ([\C], W_0, K)}$
over $\tilde U_\delta ([\C], W_0, K)$. 
For the reader's convenience, we recall briefly the 
definition of this lifted bundle: for any 
$\tilde \C  = (\tilde f, \tilde \Sigma ;\{\tilde x_i\}, \{\tilde z_j\})$
in $\tilde U_\delta ([\C], W_0, K)$, the fiber of
$E|_{\tilde U_\delta ([\C], W_0, K)}$ at $\tilde \C$
consists of all $C^{\ell - 1}$-smooth, $\tilde f^*TV$-valued $(0,1)$-forms on
$\tilde \Sigma$.
When one passes from 
$\tilde U_\delta ([\C], W_0, K)$ to another local uniformization
$\tilde U_{\delta'} ([\C'],W_0',K')$, there is an obvious bundle transition
map, which lifts the identity map on $E$, from
$E|_{\tilde U_\delta ([\C], W_0, K)}$ into 
$E|_{\tilde U_{\delta'} ([\C'], W_0', K')}$. 
Moreover, those transition maps satisfy all properties listed in section 1.
Therefore, we have a topological orbifold bundle $E$ over
$\overline {\F}^\ell _A (V, g,k)$, which is locally described by
those $E|_{\tilde U_\delta ([\C], W_0, K)}$.

The Cauchy-Riemann operator
$\Phi$ can be canonically lifted to each
local uniformization $\tilde U_\delta ([\C], W_0, K)$.

Now let us check that $\Phi: \overline {\F}^\ell _A (V, g,k)\mapsto E$
satisfy all properties (1) - (4) in the definition of
generalized Fredholm orbifold bundles.

All those $U_\delta([\C], W_0, K)$ cover
the moduli space $\Phi^{-1}(0)$ in $\overline {\F}^\ell _A (V, g,k)$.
By the Gromov Compactness Theorem (cf. [Gr], [PW], [Ye], and also 
[RT1], Proposition 3.1), $\Phi^{-1} (0)$ is compact in
$\overline {\F}_A^\ell (V,g,k)$ ($\ell \ge 2$).

For any $S \subset
\Gamma^{0,1}(\U, TV)$ with properties stated in Proposition 3.4,
we can define a bundle $E_S$ of finite rank as before, where
$W_0$, $\delta$ are small and $K$ is big. Moreover, we
assume that $S$ is invariant under the action of ${\rm Aut}(\C)$.
By Proposition 3.4, $(E_S, \Phi^{-1}(E_S))$ is a smooth
approximation of $\tilde U_\delta ([\C], W_0, K)$.
Furthermore, $(E_S, \Phi^{-1}(E_S))$ 
is invariant under the action of ${\rm Aut}(\C)$.
We denote such a smooth approximation by 
$$(\tilde E_{\delta, S} ([\C], W_0, K), \tilde 
X_{\delta, S} ([\C], W_0, K)).$$
One can easily show that all the smooth approximations of the form
$$(\tilde E_{\delta, S} ([\C], W_0, K), 
\tilde X_{\delta, S} ([\C], W_0, K))$$
are compatible with above transition maps between local uniformizations
$\{\tilde U_\delta([\C], W,K)\}$. 
Therefore, $\Phi:\overline {\F}^\ell _A (V, g,k) \mapsto E$ is weakly
smooth. Its index can be computed by the Atiyah-Singer Index
Theorem and is equal to $2c_1(V)(A) + (2n-3)(1-g) + 2k$.

We put 
$$\begin{array}{rl}
&E_{\delta, S} ([\C], W_0, K)= 
\tilde  E_{\delta, S} ([\C], W_0, K) / \Gamma,\\
&X_{\delta, S} ([\C], W_0, K) =
\tilde X_{\delta, S} ([\C], W_0, K)/ \Gamma ,\\
\end{array}
$$
where $\Gamma = {\rm Aut}(\C)$.

We claim that $\Phi^{-1}(0)$ can be covered by finitely many
smooth approximations of the form $X_{\delta, S} ([\C], W_0, K)$.
This is the same as saying
that for each $[\C]$ in $\Phi ^{-1}(0)$, there is a small
neighborhood $U$ such that for some smooth approximation
$X_{\delta, S} ([\C], W_0, K)$, 
$$[\C ]\in U\cap \Phi^{-1}(0) \subset
X_{\delta, S} ([\C], W_0, K).$$
This follows from our construction of
$X_{\delta, S} ([\C], W_0, K)$ and the following lemma.

\begin{lemma}
Let $\{f_i\}$ be a sequence of $J$-holomorphic maps with fixed
homology class $A$, then by taking a subsequence if necessary,
we may have that $f_i$ converges to some
holomorphic map $f_\infty$, which may be reducible, 
such that $||f_i||_{1,p}$ is uniformly bounded
for any $p > 2$.

\end{lemma}

\vskip 0.1in
\noindent
{\bf Proof:} This lemma was in fact essentially proved in [RT1],
section 6. It is also true for any sequence of harmonic maps (cf. [CT]).

By the Gromov Compactness Theorem, we may assume that
$f_i$ converges to $f_\infty$ in the topology of
$\overline \F_A^\ell(V,g,k)$.
Then we suffice to show that $|| f_i||_{1,p}$ is uniformly
bounded.

Let $\Sigma _i$ be the domain of $f_i$ and $q$ be an node
of $\Sigma_\infty$, which is the domain of $f_\infty$. 
Near $q$, $\Sigma _i$
can be locally described by coordinates $w_1, w_2$ with
$w_1 w_2 = \epsilon _i$ and $|w_1|, |w_2| \le 1$. Note 
that $\lim_{i\to \infty} \epsilon _i = 0$.

We may assume that $\epsilon_i > 0$ and 
$f_i(w_1, w_2)$ is very close to $q$.
Write $w_1 = s e^{\sqrt{-1} \theta}$, where $\epsilon _i \le s \le 1$.
Then the Cauchy-Riemann equation becomes
$$\frac{\partial f_i }{\partial s} + \frac{1}{s}
J(f_i) \frac{\partial f_i}{\partial 
\theta} = 0.$$
By the same arguments as in the proof of Lemma 3.9, we
can deduce that for any $p > 2$,
$$\int _{\epsilon _i 
\le s \le 1} (|w_1|^2+ |w_2|^2)^{\frac{p-2}{p}}
|\nabla f_i|^2 s ds\wedge d\theta 
\le c_p ,$$
where $c_p$ is a constant depending only on $p$.
It follows that $|| f_i||_{1,p}$ is uniformly bounded.
The lemma is proved.
\vskip 0.1in

Now let us construct a resolution $\{F_i, \psi_i\}$ of $\Phi^{-1}(0)$.
We cover $\Phi^{-1}(0)$ by smooth
approximations $\{(E_{\delta_i, S_i}([\C_i], W_{0i}, K_i),
X_{\delta_i, S_i}([\C_i], W_{0i}, K_i))\}$, where $1\le i \le m$.
For each $i$, we have 
$$X_{\delta_i, S_i}([\C_i], W_{0i}, K_i)=
\tilde X_{\delta_i, S_i}([\C_i], W_{0i}, K_i)/ \Gamma_i,$$
where $\Gamma_i$ is the automorphism group of $\C_i$ and
$$\tilde X_{\delta_i, S_i}([\C_i], W_{0i}, K_i) =
\Phi^{-1}(E_{S_i}) \subset \tilde U_{\delta _i}([\C_i], W_{0i}, K_i)=
\Map _{\delta_i}(W_i,K_i).$$

For each $i$, choose $W_{0i}'\subset W_{0i}$ such that if $W_i'$ denotes
the corresponding subset in $W_i$, then all 
the $X_{\delta_i, S_i}([\C_i], W_{0i}', K_i)$
still cover $\Phi^{-1}(0)$.

As we have seen before, any vector $v \in S_i$ 
induces a section, denoted
by $v_s$, of $E_{S_i}$ over $\Map_{\infty}(W_i)$.

Let us construct $F_i, \psi_i$ inductively. For $i=1$, we
simply define $F_1=S_1$ and 
$$\begin{array}{rl}
\psi_1: &\Map _{\infty }(W_1)\times F_1 \mapsto E_{1},\\
&\psi_1 (\tilde \C, v) = \eta _1(\Sigma_{\tilde \C}) v_s(\tilde \C).\\
\end{array}$$
Here $\eta_i$ is a smooth function on $\tilde W_i$ satisfying:
$\eta_i \equiv 1$ on $\tilde W_i'$ and $\eta_i = 0$ near $\partial \tilde W_i$.

Suppose that we have defined $F_i, \psi_i$ for $1\le i \le l-1$.
Let us define $F_{l}, \psi_{l}$. For each $i < l$ and $v\in F_i$,
$\psi_i(v)$ induces a section, say $v_{s,l,i}$, of $\tilde E_l$ over 
$\tilde U_{\delta_{l}}(\C_{l}, W_{0 l},K_{l})$. 
Let $F_l$ be the vector space spanned by
$S_l$ and $\sigma ^*(v_{s,l,i})$ ($\sigma \in \Gamma_l$,
$i< l$). We define $\psi_l(\tilde \C, v)$ to be $\eta _l v_s(\tilde \C)$ if
$v\in S_l$ and $v_{s,l,i}(\sigma(\tilde \C))$ if
$v=\sigma ^*(v_{s,l,i})$. Clearly, $\psi_l$ is $\Gamma_l$-equivariant.
Thus we can construct $\{F_i,\psi_i\}_{1\le i\le m}$.

The smooth structure of $\Phi: \overline \F_A^\ell(V,g,k)\mapsto E$
is given by all those smooth approximations
$(\tilde E_{\delta, S} ([\C], W_0, K), 
\tilde X_{\delta, S} ([\C], W_0, K))$ satisfying:
if $\U$ (resp. $\U_i$) is the universal family of curves over $\tilde W$
(resp. $\tilde W_i$)
corresponding to $W_0$ (resp. $W_{0i}$), 
then $S|_{\U\cap \U_i}$ contains
the image of $\psi_{i}$ for each $i$.
One can show that with all these 
$(\tilde E_{\delta, S} ([\C], W_0, K), 
\tilde X_{\delta, S} ([\C], W_0, K))$, 
$\{F_i, \psi_i\}$ satisfies all properties
required for a smooth resolution.
Therefore, $\Phi: \overline \F_A^\ell(V,g,k)\mapsto E$
is a generalized Fredholm orbifold bundle.

Finally, let us give the natural orientation of $\det (\Phi )$ (cf. [R], [RT1]).
We first notice that for any $\C$ representing a point
in $\Phi^{-1}(0)$,
$\det (\Phi)|_\C$ can be naturally identified with
the determinant of the linear Fredholm operator
$L_\C\Phi$, where $L_\C\Phi$ denotes
the linearization of $\Phi$ at $\C$.
It can be written as $\overline \partial _\C + B_\C$,
where $\overline \partial _\C$ is a $J$-linear operator and
$B_\C$ is an operator of zero order. It follows that
$L_\C\Phi$ can connected to $\overline \partial _\C$ through the 
canonical path $\{\overline \partial _\C + t B_\C\}_{0\le t\le 1}$,
so $\det (\Phi )$ is naturally isomorphic to the determinant
$\det (\overline \partial )$ of the family of operators 
$\{\overline \partial _\C\}_\C$. 
However, since each $\overline \partial _\C$ is J-invariant, there
is a natural orientation on $\det (\overline \partial )$. It follows that
$\det (\Phi )$ can be naturally oriented.

Proposition 2.2 is proved.

\begin{rem}

In fact, in this concrete case, we can construct the Euler 
class $e([\Phi: \overline \F_A^\ell (V,g,k) \mapsto E])$
without using Theorem 1.2. We can use the arguments in the 
proof of Theorem 1.2 and smooth approximations 
$$(\tilde E_{\delta, S} ([\C], W_0, K), 
\tilde X_{\delta, S} ([\C], W_0, K))$$
to construct a $\QQ $-cycle. This $\QQ$-cycle will lie
in a finite covering of $\Phi^{-1}(0)$ by finitely
dimensional smooth approximations. Here we do need
smooth properties of $(\tilde E_{\delta, S} ([\C], W_0, K),
\tilde X_{\delta, S} ([\C], W_0, K))$
during changes of local uniformizations.  

\end{rem}

\vskip 0.1in
\noindent
{\bf Proof of Proposition 2.3:} The proof is identical to that of 
Proposition 2.2. So we omit the details.

Let $\Phi: \overline F_A^\ell (V,g,k)  \mapsto E$ be the
generalized Fredholm orbifold bundle as above,
and $\Phi': \overline F_A^\ell (V,g,k)  \mapsto E$ is another one
induced by the almost complex structure $J'$.

Let $\{J_t\}$ be the family of almost complex structures joining $J$ to $J'$.
Consider 
$$\begin{array}{rl}
\Psi: [0,1] \times \overline F_A^\ell (V,g,k) & \mapsto E,\\
(t, (f, \Sigma; \{x_i\})) &\mapsto df + J_t(f)\cdot df \cdot j.\\
\end{array}
$$
Then $\Psi |_{\{0\}\times \overline F_A^\ell (V,g,k) } = \Phi$
and $\Psi |_{\{1\}\times \overline F_A^\ell (V,g,k) } = \Phi'$.

Using the same arguments as above, one can prove that
$\Psi: [0,1]\times \overline F_A^\ell (V,g,k)  \mapsto E$
is a generalized Fredholm orbifold bundle.
Moreover, one can equip this bundle a weakly smooth
structure which restricts to
the given smooth structures of $\Phi$ and $\Phi'$
along the boundary $\{0, 1\}\times \overline F_A^\ell (V,g,k) $.

This follows that $\Phi$ is homotopic to $\Phi'$. That is 
just what Proposition 2.3 claims.

\vskip 0.1in

\section{One more example}
\label{sec:4}

In this section, we consider a simpler example:
the Seiberg-Witten invariants of
4-manifolds. The Seiberg-Witten invariants
have found many striking applications in the study
of 4-dimensional topology (cf. [Wi], [Ta], [KM], [KST]).
Here we just give a different approach to defining the Seiberg-Witten
invariants, which seems to be of independent interest.

We first fix the notation we will use.
Let $X$ be a compact oriented smooth 4-manifold
and let $c$ be a $\spin^c$ structure on $X$ with the associated
$\spin^c$ bundles $\wp$ and $\wm$. Let
\begin{displaymath}
  \rho: \Lamp\otimes\CC\lra\sl(\wp)
\end{displaymath}
be the isomorphism induced by the Clifford multiplication,
where $\sl(\wp)$ is the associated $PSL(2,\CC)$ bundle of $\wp$.
There is also a pairing
\begin{displaymath}
  \wp\times\bar{W}^+\lra\sl(\wp)
\end{displaymath}
that is modeled on the map $\CC\times\bar\CC^2\to\sl(\CC^2)$ sending
$(v,w)$ to $i(v\bar w^t)_0$, where the subscript means the traceless
part.

Now the Seiberg-Witten invariants is defined as follows.
We first fix a Riemannian metric $g$.
Let $L$ be the determinant line bundle of $\wp$ and let $\A (\L)$ be the
space of unitary connections on $L$. Then $A\in\A (\L)$ induces a
Dirac operator $\Gamma(\wp)\to\Gamma(\wm)$. Now let $\tilde\B$ be
the Banach manifold $\A(L)\times\Gamma(\wp)$ and let $\tilde\E$
be the constant vector bundle over $\tilde\B$ with fiber
$\Gamma(\wp)\times\Gamma(\sl(\wp))$. We define a section
$\tilde f: \tilde\B\to\tilde\E$ via
\begin{displaymath}
  \tilde f(\varphi, A)=(D_A\varphi, \rho( F_A^+)-i\sigma(\varphi,\varphi)).
\end{displaymath}
Now let $p_0\in X$ be  fixed and let $\G_0=\Map_{p_0}(X,S^1)$
be the pointed gauge group of $L$. Note that $\G_0$ also acts on
$\Gamma(\wp)$ via scalar automorphism of $\wp$. Hence $\G_0$ acts freely on
$\tilde B$  and it lifts to an action on $\tilde \E$. We let
$\B=\tilde\B/\G_0$ and $\E=\tilde\E/\G_0$. Since $\tilde f$ is
equivariant under $\G_0$, $\tilde f$ descends to a section
\begin{displaymath}
  f: \E\lra\B\,.
\end{displaymath}
Note that $f^{-1}(0)$ is compact.
Now let $\G$ be the full gauge group. Then $\G/\G_0\cong S^1$
acts on $\E\to\B$ and the section $f$ is $S^1$-equivariant as well. The
Seiberg-Witten invariants of $X$ is the $S^1$-equivariant version of
Euler class of $[f:\B\to\E]$  defined in section 1. More precisely,
The universal line bundle on $\A(L)$ descends to a complex
line bundle $\L$ on $\B$ and the Seiberg-Witten invariants of
$X$ is
\begin{displaymath}
  SW: H^2(X,\ZZ)\lra \ZZ
\end{displaymath}
defined by
\begin{displaymath}
  SW(L)=<\Euler(\E,f)^{S^1}, c_1(\L)^k>\,,
\end{displaymath}
where
\begin{displaymath}
  k={\frac{1}{4}}(c_1(L)^2-(2\chi+3\sigma)),
\end{displaymath}
is the Fredholm index of $f/S^1: \B/S^1\to \E/S^1$.

In case the zero locus $f^{-1}(0)$ is disjoint from the fixed point
set of $S^1$, which is
\begin{displaymath}
  \B\uso=\A(L)/\G_0\times\{0\}\sub \A(L)\times\Gamma(\wp)/\G_0,
\end{displaymath}
then we can work with $[f/S^1:\B/S^1\to\E/S^1]$. Let $\B'=(\B-\B^{S^1})/S^1$,
let $\E'=(\E|_{\B'})/S^1$ and let $f'=(f|_{\B'})/S^1$. Then $[f': \B'\to \E']$
is a Fredholm operator as defined in section 1. The
Seiberg-Witten
invariant then is
$$SW(L)=<e[f': \B'\to \E'], c_1(\L')^k>
$$
where $\L'$ is the descend of $\L|_{\B-\B^{S^1}}$ to $\B'$.

We now look at the general case. Let
\begin{displaymath}
  \E|_{\bso}=\oplus_{i=\-\infty}^{\infty}\F_i
\end{displaymath}
be the spectral decomposition of the restriction of $\E$ to $\bso$.
Namely, $\F_i\sub\E|_{\bso}$ is
$S^1$-invariant and the $\so$ action on $\F_i$ has weight $i$.
Then $f|_{\bso}$ factor through $\F_0\sub\E|_{\bso}$. We denote this
section by $f_0$. Let $k+1$ be the Fredholm index of $f$ and let
$l$ be the Fredholm index of $f_0: T_z\bso \to \F_{0,z}$.

\begin{lemma}
  Assume $l<0$, then any $\so$-equivariant Fredholm section
$f:\B\to\E$ is homotopic to an $\so$-equivariant Fredholm section
$g:\B\to\E$ so that $g^{-1}(0)\cap\B^{S^1}=\emptyset$.
\end{lemma}

\noindent
{\bf Proof:}
The proof is straightforward. We first look at the
the restriction of $f$ to $\B^{S^1}$. As we mentioned, it factor
through $\F_0$. Let $h: \B^{S^1}\to \F_0$ be this map. Then since
$dh$ has negative Fredholm index, by Theorem 1.1, $h$ is homotopic to
$\tilde h: \B^{S^1}\to\F_0$ so that its vanishing locus is empty.
Clearly, for some $S^1$ invariant neighborhood $U$ of $\B^{S^1}\sub
\B$, we can extend this homotopy, and thus $\tilde h$, within
the category of Fredholm operators, to an $S^1$-equivariant
$g:\B \to\E$ so that the restriction of $g$
to $\B^{S^1}$ is $\tilde h$ and $g|_{\B-U}=f|_{\B-U}$.
This proves the Lemma.

After having $g$ given by the Lemma, we reduce the situation to when
$g^{-1}(0)\cap  \B^{S^1}=\emptyset$. Thus as before, we can define
the equivariant Euler class $e[f: \B^{S^1}\to\E]^{S^1}$
represented by a smooth $k$ dimensional submanifold in $(\B-\B^{S^1})/S^1$
to be the class $e[g/S^1: \B'\to \E']$. When $l<-1$, then the
above argument shows that any two such representatives of
$e[f: \B\to\E]^{S^1}$ in $\B'$ are coborbant to each other in $\B'$.
Therefore, they represent a well-defined corbordism class.

We now apply the above construction to the Seiberg-Witten invariant, the
fixed point
set $\bso$ in $\B$ is $\A(L)\times\{0\}$. The $\F_0\sub\E|_{\bso}$
in this case is the subbundle $\Gamma(\sl(\wp))$ and the restriction section
is $\rho(F_A^+)$, whose Fredholm index is $-b_2^+$. Therefore, when $b_2^+>1$,
the Seigerg-Witten invariant is well defined and
can be represented by a smooth submanifold in $(\B-\bso)/S^1$.


\begin{thebibliography}{L3}

\bibitem[B]{B} K.~Behrend, Gromov-Witten invariants in algebraic geometry, preprint, 
1996.

\bibitem[BF]{BF} K.~Behrend and B.~Fantechi, The intrinsic normal cone,
preprint, 1996 

\bibitem[CT]{CT} J.Y.~Chen and G.~Tian,
Compactification of moduli space of harmonic mappings, preprint 1996.

\bibitem[FO]{FO} K. ~Fukaya and K.~ Ono, 
Arnold conjecture and Gromov-Witten invariants,
preprint, 1996.

\bibitem[Gr]{Gr} M. ~Gromov, Pseudo holomorphic curves in symplectic manifolds,
                Invent. math., 82 (1985), 307-347.

\bibitem[KM]{KM} M.~ Kontsevich and Y. ~Manin, GW classes, Quantum
cohomology and enumerative geometry,
Comm. Math. Phys. vol. 164, 1994, 525-562.

\bibitem[KST]{KST} J.~Morgan, Szabo and C. Taubes, preprint.

\bibitem[LT1]{LT1} J. ~Li and G. ~Tian,
Virtual moduli cycle and Gromov-Witten invariants of algebraic varieties,
preprint, 1996.

\bibitem[LT2]{LT2} J. ~Li and G. ~Tian,
Algebraic and symplectic geometry of virtual moduli cycles, submitted to
Proc. of AMS summer school, 1995, Santa Cruz.

\bibitem[LiuT]{LiuT} G. ~Liu and G.~ Tian,
Arnold conjecture for general symplectic manifolds, in preparation.

\bibitem[Mu]{Mu} D.~ Mumford, 
Towards an enumerative geometry of the moduli space of curves,
Arithmetic and Geometry II, Progress in Mathematics 36,
1983, 271-326.

\bibitem[PW]{PW}  T. ~Parker and J. ~Woflson, A compactness theorem for
                Gromov's  moduli space, J. Geom. Analysis, 3 (1993).

\bibitem[Ru]{Ru} Y. ~Ruan, Topological Sigma model and Donaldson type 
invariants in Gromov theory, to appear in Duke J. Math.

\bibitem[RT1]{RT1} Y. ~Ruan and G. ~Tian,  
A mathematical theory of quantum cohomology, J. Diff. Geom.. vol 42, no. 2,
1995.

\bibitem[RT2]{RT2} Y. ~Ruan and G. ~Tian, Higher genus symplectic 
invariants and sigma model coupled with gravity, preprint, 1995.

\bibitem[Si]{Si} B. ~Siebert, 
Gromov-Witten invariants for general symplectic manifolds,
preprint, 1996.

\bibitem[Ta]{Ta} C. ~Taubes, The Seiberg-Witten and the Gromov invariants,  
                 preprint, 1995.

\bibitem[T]{T} G. ~Tian, Quantum cohomology and its associativity,
Proc. of 1st Current developments in Mathematics, Cambridge, 1995.

\bibitem[Wi]{Wi}  E. ~Witten, Monopoles and 4-manifolds, Math. Res. Lett., vol. 1,
1994, 769-796.

\bibitem[Wi1]{Wi1}  E. ~Witten,  Topological sigma models, 
Comm. Math. Phys.. vol. 118, 1988.

\bibitem[Wi2]{Wi2} E. ~Witten,  Two dimensional gravity and intersection
theory on moduli space,  Surveys in Diff. Geom..
vol  1, 1991, 243-310.

\bibitem[Ye]{Ye} R. ~Ye, Gromov's compactness theorem for pseudo-holomorphic
curves, Trans. Amer. math. Soc., 1994.


\end{thebibliography}
\end{document}